%% file: main.tex
\documentclass[onefignum,onetabnum]{siamart171218}
\usepackage{amsmath}
\usepackage{lineno}
\input{ex_shared}
\usepackage{graphicx}

\usepackage{tikz}
\usetikzlibrary{shapes.geometric}
\usetikzlibrary{arrows.meta,positioning,fit,calc}

\usepackage{subcaption}
\usepackage{booktabs}

\usepackage{soul}
\usepackage{comment}

\newcommand{\rev}[1]{\textcolor{black}{#1}}

\ifpdf
\hypersetup{
  pdftitle={A Hybrid Quantum-Classical Coarse-to-Fine Finite Element Solver for Heterogeneous Darcy Flow},
  pdfauthor={Kapil Chawla, Sanghyun Lee, and Jonathan Valyou}
}
\fi

\begin{document}

\maketitle

\begin{abstract}
Quantum linear system algorithms offer
theoretical potential for solving large sparse systems from
PDE discretizations. However, applying them directly to fine-scale problems
remains impractical on near-term hardware because of limited qubits,
restricted circuit depth, and noise.
\rev{They also prepare a quantum state rather than solution coefficients,
whose recovery requires measurement and reconstruction.}

{\color{black}In this work, we develop a hybrid quantum-classical coarse-to-fine finite
element framework for heterogeneous Darcy flow. We investigate whether a
measured quantum coarse solution improves fine-grid initialization and reduces
subsequent classical solver work. {\color{black}Rather than solve the full fine-scale system
on a quantum device, we apply a method inspired by the Variational Quantum
Linear Solver (VQLS) to the Galerkin coarse system.} {\color{black}We reconstruct the
approximate coarse solution from measurements, prolongate it to the fine grid,
and use it to initialize a classical preconditioned iterative solve.} Numerical experiments on steady heterogeneous Darcy flow and a
time-dependent SPE10-based problem show that measured quantum coarse
solutions can improve fine-grid initialization, while classical correction
recovers accurate pressure fields. The effect on classical solver iterations
depends on the problem. All quantum computations were simulated using Qiskit Aer.}
\end{abstract}



\section{Introduction}

Darcy flow models fluid motion in porous media, with applications in
groundwater flow, geothermal energy extraction, and carbon
sequestration. Geological permeability can vary by orders of
magnitude or be discontinuous
\cite{EfendievHou2009,Durlofsky1991}. Resolving the resulting local
flow features requires fine meshes, and the heterogeneity can produce
poorly conditioned linear systems~\cite{Saad2003,EfendievHou2009}.
A coarse discretization may capture the broad pressure trend but miss
effects of fine-scale permeability \cite{ChawlaLeeShin2026PAM}.
{\color{black}The steady Darcy problem is elliptic for positive permeability;
pressure storage adds a time derivative, giving a parabolic
extension~\cite{ChenHuanMa2006}. The steady and implicit
time-step finite element systems are real symmetric positive
definite (SPD).}

Quantum linear system algorithms offer theoretical potential for
operations on large sparse systems
\cite{Balducci2022,Schalkers2024,Morales2024QLSsurvey}. Under suitable
assumptions, the Harrow-Hassidim-Lloyd algorithm~\cite{Harrow2009}
prepares a state encoding a linear-system solution with
polylogarithmic dependence on the system dimension.
{\color{black}Later work improved the dependence on precision~\cite{Childs2017}}, and quantum singular value transformation provides a broader
framework for quantum linear algebra~\cite{Gilyen2019}. The value of
these methods depends on conditioning, access to the matrix and
right-hand side, and the information needed from the solution
\cite{Aaronson2015}. {\color{black}In particular, a quantum solution state does not
directly provide the classical pressure coefficients; recovering them
requires measurement and reconstruction~\cite{Henderson2024}.}

Current devices also limit the number of qubits and circuit depth and
are affected by gate and readout noise~\cite{Preskill2018}.
Variational methods such as the Variational Quantum Linear Solver
(VQLS)~\cite{BravoPrieto2023} approximate a solution with a
parameterized circuit and objectives estimated from measurements.
{\color{black}These constraints motivate placing the quantum solve on a small
coarse system while retaining a classical solve on the fine grid
\cite{Pollachini2021,Jacob2025}. 
The question is whether its solution remains useful for the fine-grid
solve after measurement and classical reconstruction.}
{\color{black} Therefore, our motivation is to examine the algorithmic role of an approximate quantum coarse solve within a hybrid coarse-to-fine PDE workflow.}

{\color{black}Quantum methods have been studied for discretized PDEs
\cite{Lubasch2020,Liu2021,Sato2021,Fang2023,Deiml2024,Arora2025}.
For subsurface flow, Golden et al.~\cite{Golden2022} propose a
quantum-compatible inverse-Laplacian preconditioner. Henderson
et al.~\cite{Henderson2023} study quantum fracture-network solvers,
including a variational method; their later work~\cite{Henderson2024}
addresses input preparation and extraction of solution information. Here, we focus
on a measured quantum coarse solution as an initial approximation
for a classical fine-grid Darcy solver.}

{\color{black}In this work, we form the coarse system by Galerkin projection of the fine-grid
finite element system and pad it to dimension $2^n$ when needed for
an $n$-qubit register. A VQLS-inspired circuit approximates the
coarse solution direction using an energy-based objective estimated
from measurements~\cite{BravoPrieto2023,Sato2021}. Tomography and
classical postprocessing reconstruct and rescale the coarse
coefficient vector. We then prolongate this vector to the fine grid,
add the boundary lifting, and use the result to initialize an
algebraic-multigrid-preconditioned conjugate-gradient (AMG-PCG)
solve~\cite{Saad2003,LiZikatanovZuo2025}.}

{\color{black}Beyond constructing the solver, we examine what the measured quantum
coarse solution contributes to the fine-grid solve:
\begin{itemize}
\item We use Galerkin orthogonality~\cite{Briggs2000} to separate the
initial fine-grid energy error into coarse-space approximation error
and error introduced by the quantum coarse solve and reconstruction
(see~\eqref{eq:initial_energy_decomposition} in
Section~\ref{sec:fine_scale_corr}).

\item We compare measured quantum, zero, and exact Galerkin
initializations using initial $L^2$ and energy errors and subsequent
AMG-PCG iterations. This comparison identifies whether a better
initial approximation changes the classical correction.

\item We study steady and time-dependent heterogeneous Darcy flow
under finite sampling and modeled noise, accounting for the
preparation and measurements needed to recover the coarse solution.
\end{itemize}
}

Section~\ref{sec:model_problem} presents the finite element and
Galerkin formulations. Section~\ref{sec:quantum_solver} describes the
quantum coarse solve and classical reconstruction.
Section~\ref{sec:fine_scale_corr} presents fine-grid correction.
Section~\ref{sec:numerical_experiments} reports the numerical
experiments, and Section~\ref{sec:conclusion} concludes the paper.

\section{Model Problem and Finite Element Discretization}
\label{sec:model_problem}

We define the steady Darcy problem and its fine-grid finite element
system, then construct the Galerkin coarse system used by the quantum
solver. The time-dependent discretization is specified with
Experiment \ref{sec:ex2}.

Let $\Omega\subset\mathbb{R}^2$ be a bounded Lipschitz domain. Its
boundary is divided into disjoint Dirichlet and Neumann parts,
$\partial\Omega=\Gamma_D\cup\Gamma_N$, with $|\Gamma_D|>0$.
The steady pressure $p:\Omega\to\mathbb{R}$ satisfies
\begin{subequations}\label{eq:model}
\begin{alignat}{2}
-\nabla\cdot\bigl(K(\mathbf{x})\nabla p\bigr)&=f
    &&\quad\text{in }\Omega,\\
p&=g_D&&\quad\text{on }\Gamma_D,\\
\bigl(K(\mathbf{x})\nabla p\bigr)\cdot\mathbf{n}&=g_N
    &&\quad\text{on }\Gamma_N.
\end{alignat}
\end{subequations}
Here $\mathbf{x}$ is the spatial position, $\mathbf{n}$ is the outward
unit normal, $f$ is the source term, $g_D$ is the prescribed pressure,
and $g_N$ is the prescribed normal component of $K\nabla p$.
{\color{black}We use the isotropic permeability tensor
$K(\mathbf{x})=\kappa(\mathbf{x})I$, where $I$ is the $2\times2$
identity matrix. The scalar permeability $\kappa$ may be discontinuous,
but we assume
\[
0<\kappa_{\min}\leq\kappa(\mathbf{x})\leq\kappa_{\max}<\infty
\qquad\text{for almost every }\mathbf{x}\in\Omega.
\]}
The ratio $\kappa_{\max}/\kappa_{\min}$ measures permeability
contrast. Large contrast can make the discrete system difficult to
solve.

\subsection{Weak formulation and fine-grid system}

Let $H^1(\Omega)$ be the Sobolev space of square-integrable functions
with square-integrable first derivatives. Boundary equalities below
are understood in the trace sense. Define the trial and test spaces
\[
V_g=\{v\in H^1(\Omega):v=g_D\text{ on }\Gamma_D\},\qquad
V_0=\{v\in H^1(\Omega):v=0\text{ on }\Gamma_D\}.
\]
Multiplying \eqref{eq:model} by a test function and integrating by
parts gives: find $p\in V_g$ such that
\begin{equation}\label{eq:darcy_weak}
a(p,v)=\ell(v)\qquad\text{for every }v\in V_0,
\end{equation}
where
\[
a(u,v)=\int_\Omega K\nabla u\cdot\nabla v\,d\mathbf{x},\qquad
\ell(v)=\int_\Omega fv\,d\mathbf{x}
          +\int_{\Gamma_N}g_Nv\,ds.
\]
The sign of the Neumann term follows the boundary condition in
\eqref{eq:model}.

Let $\mathcal{T}_h$ be a conforming triangulation of $\Omega$ with
mesh size
$h:=\max_{T\in\mathcal{T}_h}\operatorname{diam}(T)$, where
$\operatorname{diam}(T)$ is the diameter of triangle $T$. We use
continuous piecewise linear finite elements,
\[
V_h=\{v_h\in H^1(\Omega):v_h|_T\in\mathbb{P}_1(T)
       \text{ for every }T\in\mathcal{T}_h\},
\]
where $\mathbb{P}_1(T)$ is the space of affine functions on a
triangle $T$. Let $V_{h,g}=V_h\cap V_g$ and $V_{h,0}=V_h\cap V_0$;
the prescribed Dirichlet values are imposed at boundary nodes.
{\color{black}We write $K_h=\kappa_hI$ for the fine-mesh permeability used in
assembly and assume $\kappa_h$ is uniformly positive and bounded. The fine
fields used in Experiments \ref{sec:ex1} and \ref{sec:ex2} are shown in
Figure~\ref{fig:pressure_comparison}(a) and
Figure~\ref{fig:ex2_permeability}, respectively. }

With $a_h(u,v)=\int_\Omega K_h\nabla u\cdot\nabla v\,d\mathbf{x}$,
the discrete problem is: find $p_h\in V_{h,g}$ such that
\begin{equation}\label{eq:darcy_discrete_weak}
a_h(p_h,v_h)=\ell(v_h)\qquad\text{for every }v_h\in V_{h,0}.
\end{equation}
Let $N_h^{\mathrm{tot}}=\dim V_h$ be the total number of fine-grid
nodal degrees of freedom, including boundary nodes, and let
$N_h^D$ count those on $\Gamma_D$. The remaining
$N_h=N_h^{\mathrm{tot}}-N_h^D=\dim V_{h,0}$ degrees of freedom are
unknown. Let $\{\phi_i\}_{i=1}^{N_h}$ be the nodal basis associated
with these free nodes. Eliminating the prescribed Dirichlet values
gives
\begin{equation}\label{eq:linear_system}
A_h\mathbf{p}_h=\mathbf{b}_h,
\end{equation}
where $\mathbf{p}_h\in\mathbb{R}^{N_h}$ contains the unknown
pressures at free nodes, $\mathbf{b}_h\in\mathbb{R}^{N_h}$ is the
reduced load vector, and
\[
A_h\in\mathbb{R}^{N_h\times N_h},\qquad
(A_h)_{ij}=a_h(\phi_j,\phi_i).
\]
The Dirichlet values are prescribed separately; they are not entries
of $\mathbf{p}_h$. The reduced load includes the source and Neumann
terms and the contribution from eliminating those prescribed values.
The positivity of $\kappa_h$ and $|\Gamma_D|>0$ make $A_h$ symmetric positive
definite (SPD).

\subsection{Galerkin coarse system}\label{sec:model_problem_coarse}

{\color{black}Let $\mathcal{T}_H$ be a coarse triangulation with mesh size
$H:=\max_{E\in\mathcal{T}_H}\operatorname{diam}(E)>h$, and assume
$\mathcal{T}_h$ is its conforming refinement. Define the
continuous piecewise linear coarse space $V_H$ in the same way as
$V_h$. Its total nodal count is $N_H^{\mathrm{tot}}=\dim V_H$;
after excluding prescribed Dirichlet nodes, the free coarse count is
$N_H=\dim V_{H,0}$, where
$V_{H,0}=V_H\cap V_0\subset V_{h,0}$. Thus $N_H<N_h$.

Choose a boundary lifting $G_h\in V_{h,g}$ that matches $g_D$ on
$\Gamma_D$, and let $\mathbf{g}_h\in\mathbb{R}^{N_h}$ contain its
values only at the free fine nodes. Then $p_h-G_h\in V_{h,0}$ has
free-node coefficient vector $\mathbf{p}_h-\mathbf{g}_h$. We
approximate this homogeneous correction by $p_H\in V_{H,0}$.
For homogeneous Dirichlet data, we take $G_h=0$. The Galerkin
problem is
\begin{equation}\label{eq:coarse_weak}
a_h(p_H,v_H)=\ell(v_H)-a_h(G_h,v_H)
\qquad\text{for every }v_H\in V_{H,0}.
\end{equation}

Let $\{\Phi_j\}_{j=1}^{N_H}$ be the nodal basis of $V_{H,0}$ and write
$p_H=\sum_{j=1}^{N_H}(\mathbf{p}_H)_j\Phi_j$.
The vector $\mathbf{p}_H\in\mathbb{R}^{N_H}$ contains the
homogeneous coarse coefficients. Because the spaces are nested,
each $\Phi_j$ can be expressed in the fine basis. This defines the
prolongation matrix $P\in\mathbb{R}^{N_h\times N_H}$ by
\[
\Phi_j=\sum_{i=1}^{N_h}P_{ij}\phi_i,
\qquad j=1,\ldots,N_H.
\]
The matrix $P$ maps coarse coefficients to values at the free fine
nodes. Galerkin projection of \eqref{eq:linear_system} gives
\begin{equation}\label{eq:coarse_linear_system}
A_H\mathbf{p}_H=\mathbf{b}_H,\qquad
A_H=P^TA_hP,\qquad
\mathbf{b}_H=P^T(\mathbf{b}_h-A_h\mathbf{g}_h).
\end{equation}
The subtraction removes the free-node part of the boundary lifting
before restricting the fine-grid residual. Since $A_h$ is SPD and
$P$ has full column rank, $A_H\in\mathbb{R}^{N_H\times N_H}$ is
also SPD. Moreover,
\[
(A_H)_{ij}=a_h(\Phi_j,\Phi_i)
=\int_\Omega K_h\nabla\Phi_j\cdot\nabla\Phi_i\,d\mathbf{x}.
\]
Thus, the coarse operator retains the fine-mesh permeability through
the Galerkin projection.}

Denote the exact solution of \eqref{eq:coarse_linear_system} by
$\mathbf{p}_H^\star$ and the measured, reconstructed approximation
from the quantum solver by $\mathbf{p}_H^q$. The corresponding
coarse functions lie in $V_{H,0}$. Their free-node fine-grid
approximations are $\mathbf{g}_h+P\mathbf{p}_H^\star$ and
$\mathbf{g}_h+P\mathbf{p}_H^q$, respectively. The prescribed
Dirichlet values are restored when constructing the full pressure
field.

In the time-dependent experiment, backward Euler produces an SPD
fine-grid system at each time step. The same Galerkin projection
defines its coarse system; the time-discrete operator and right-hand
side are given with the experiment in
Section~\ref{sec:ex2}.


\section{Quantum Coarse Solver}
\label{sec:quantum_solver}

We apply a variational quantum procedure to the Galerkin coarse system
\eqref{eq:coarse_linear_system}. In this section, we write it as
\begin{equation}\label{eq:coarse_system_generic}
A\mathbf u=\mathbf b,
\qquad
A=A_H,\quad \mathbf u=\mathbf p_H,\quad
\mathbf b=\mathbf b_H,\quad N=N_H.
\end{equation}
Thus $A\in\mathbb R^{N\times N}$ is symmetric positive definite
(SPD), $\mathbf b\in\mathbb R^N$ is its right-hand side, and
$\mathbf u\in\mathbb R^N$ contains the coefficients of the
homogeneous coarse pressure correction.

The circuit prepares a state representing an approximate solution
direction. Measurements and classical postprocessing recover an
explicit coarse coefficient vector, which is then used for
fine-grid initialization.

{\color{black}The procedure is inspired by the Variational Quantum Linear Solver
(VQLS)~\cite{BravoPrieto2023}, but uses an objective derived from the SPD
quadratic energy~\cite{Sato2021} rather than a standard VQLS cost. Guarantees for
the standard cost are therefore not assumed for this objective.}

\subsection{Padding and input preparation}

An $n$-qubit state has $m=2^n$ amplitudes. We choose
$n=\lceil\log_2 N\rceil$, so that $m\geq N$. If $m>N$, we pad the
coarse system as
\begin{equation}\label{eq:padded_system}
\widetilde A_\gamma=
\begin{pmatrix}A&0\\0&\gamma I_{m-N}\end{pmatrix},
\qquad
\widetilde{\mathbf b}=
\begin{pmatrix}\mathbf b\\\mathbf 0\end{pmatrix},
\end{equation}
where $I_{m-N}$ is the identity on the added components,
$\mathbf 0\in\mathbb R^{m-N}$, and
$\lambda_{\max}(A)$ is the largest eigenvalue of $A$. The choice
$\gamma=\lambda_{\max}(A)>0$ makes $\widetilde A_\gamma$ SPD
without changing the spectral condition number of $A$; the exact
padded solution is $(\mathbf u^\top,\mathbf 0^\top)^\top$.
If $m=N$, we set $\widetilde A_\gamma=A$ and
$\widetilde{\mathbf b}=\mathbf b$.

We scale the operator by its largest eigenvalue,
\begin{equation}\label{eq:scaled_quantum_operator}
s=\lambda_{\max}(\widetilde A_\gamma),
\qquad \overline A=\widetilde A_\gamma/s.
\end{equation}
The eigenvalues of $\overline A$ lie in $(0,1]$. Padding and scaling
preserve the solution direction. The eigenvalues used here are
computed classically for the small coarse matrix.

{\color{black}
The right-hand side $\mathbf b=P^T(\mathbf b_h-A_h\mathbf g_h)$
is assembled classically from the fine-grid problem. Although $A_h$ is sparse, $\mathbf b$ need not have sparse or
simply structured amplitudes: its entries depend on the forcing,
boundary lifting, fine-scale permeability, and coarse basis. Its finite
element origin does not make state preparation free, and we assume
no efficient data-loading oracle. If $\mathbf b=\mathbf 0$, the
coarse solution is zero and no quantum solve is needed. Otherwise,
define the normalized padded right-hand-side state
\[
|\widehat{\mathbf b}\rangle
=\frac{\widetilde{\mathbf b}}{\|\widetilde{\mathbf b}\|_2}.
\]
Here $|0\rangle^{\otimes n}$ denotes the joint all-zero state of
the $n$ qubits. From the known entries of $\widetilde{\mathbf b}$,
we synthesize a unitary $U_b\in\mathbb C^{m\times m}$ satisfying
\begin{equation}\label{eq:rhs_preparation}
U_b|0\rangle^{\otimes n}=|\widehat{\mathbf b}\rangle.
\end{equation}
We use Qiskit's \texttt{StatePreparation} operation. Generic known-amplitude preparation methods include isometry-based
synthesis~\cite{Iten2016} and uniformly controlled
rotations~\cite{Mottonen2005}; their gate costs can grow as
$O(2^n)$ for $n$ qubits. Its inverse
$U_b^\dagger$ is executed in each overlap-measurement circuit, so
preparation contributes to circuit depth and noise exposure.

\subsection{Parameterized trial state}

{\color{black}The ansatz prepares a normalized trial state
\begin{equation}\label{eq:trial_state}
|\psi(\boldsymbol\theta)\rangle
=U(\boldsymbol\theta)|0\rangle^{\otimes n},
\end{equation}
where $\boldsymbol\theta=\{\theta_{\ell,q}\}$ contains one rotation
angle for each layer $\ell=0,\ldots,L-1$ and qubit
$q=0,\ldots,n-1$. Thus there are $nL$ parameters. Each of the $L$
layers applies $R_y$ rotations to all qubits, followed by a chain
of CNOT gates between neighboring qubits. For an angle $\vartheta$,
\[
R_y(\vartheta)=
\begin{pmatrix}
\cos(\vartheta/2)&-\sin(\vartheta/2)\\
\sin(\vartheta/2)&\cos(\vartheta/2)
\end{pmatrix}.
\]
Because $R_y$ and CNOT have real matrix representations, the ideal
trial amplitudes are real, matching the real solution of
\eqref{eq:coarse_system_generic}. The rightmost gate acts first;
therefore, if $U_{\rm rot}^{(\ell)}$ and
$U_{\rm ent}^{(\ell)}$ denote the combined rotations and CNOTs in
layer $\ell$, respectively, then
\[
U^{(\ell)}(\boldsymbol\theta)
=U_{\rm ent}^{(\ell)}U_{\rm rot}^{(\ell)},
\qquad
U(\boldsymbol\theta)
=U^{(L-1)}\cdots U^{(0)}.
\]

\begin{figure}[!h]
\centering
\includegraphics[width=\linewidth]
{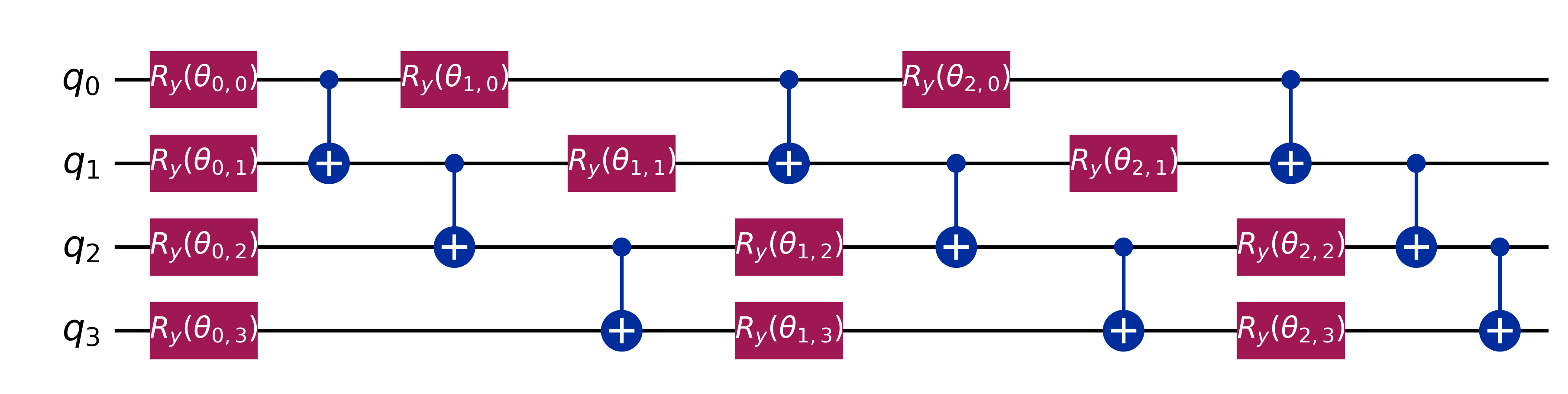}
\caption{Four-qubit trial-state circuit with three layers. Each
layer applies $R_y$ rotations followed by a nearest-neighbor CNOT
chain.}
\label{fig:quantum_ansatz}
\end{figure}

The nearest-neighbor structure follows hardware-efficient ansatz
designs~\cite{Kandala2017,BravoPrieto2023}. It uses $L(n-1)$
CNOT gates before backend compilation; state preparation and
measurement circuits add further gates. The layer count is an
experimental choice, and the number of parameters alone does not
guarantee an accurate coarse solve}.

\subsection{Energy objective and measurements}

{\color{black}For a real vector $\mathbf v\in\mathbb R^m$, the padded SPD system
has quadratic energy
\[
E(\mathbf v)=\tfrac12\mathbf v^\top\widetilde A_\gamma\mathbf v
-\widetilde{\mathbf b}^{\top}\mathbf v.
\]
Its minimizer is the padded solution. Writing
$\mathbf v=\alpha\psi(\boldsymbol\theta)$ separates the unknown
scale $\alpha\in\mathbb R$ from the normalized trial direction.
Minimizing over $\alpha$ gives
\[
\min_{\alpha\in\mathbb R}E(\alpha\psi)
=-\frac{(\widetilde{\mathbf b}^{\top}\psi)^2}
        {2\psi^\top\widetilde A_\gamma\psi}
=-\frac{\|\widetilde{\mathbf b}\|_2^2}{2s}
  \frac{p_0(\boldsymbol\theta)}{D(\boldsymbol\theta)},
\]
where $\psi$ is the real amplitude vector of
$|\psi(\boldsymbol\theta)\rangle$.} {\color{black}We minimize the energy-based objective
\begin{equation}\label{eq:quantum_objective}
\mathcal J(\boldsymbol\theta)
=-\frac{p_0(\boldsymbol\theta)}{D(\boldsymbol\theta)}
 +\eta_{\rm pad}p_{\rm pad}(\boldsymbol\theta),
\end{equation}
where $\eta_{\rm pad}\geq0$ penalizes padded amplitudes.
Unlike standard VQLS, which aligns $\overline A\psi$ with
$\widehat{\mathbf b}$~\cite{BravoPrieto2023}, this term
comes from minimizing SPD energy over scale~\cite{Sato2021}.

The three quantities in \eqref{eq:quantum_objective} are
\begin{align*}
p_0(\boldsymbol\theta)
 &=|\langle\widehat{\mathbf b}|\psi(\boldsymbol\theta)\rangle|^2,\\
D(\boldsymbol\theta)
 &=\langle\psi(\boldsymbol\theta)|\overline A
   |\psi(\boldsymbol\theta)\rangle,\\
p_{\rm pad}(\boldsymbol\theta)
 &=\langle\psi(\boldsymbol\theta)|(I_m-\Pi_N)
   |\psi(\boldsymbol\theta)\rangle,
\end{align*}
where $I_m$ is the $m\times m$ identity and
$\Pi_N=\operatorname{diag}(I_N,0_{m-N})$ projects onto the first
$N$ amplitudes; $0_{m-N}$ is the square zero matrix on the padded
components. If $m=N$, then $p_{\rm pad}=0$.
For the real normalized trial states, $p_0$ is the squared overlap
with the right-hand-side state and $D=\psi^\top\overline A\psi$.}

{\color{black}Each circuit execution, or shot, produces one measurement outcome;
$S$ denotes the shots used for a measurement setting.
The three terms are estimated as follows.

\paragraph{Overlap}
Prepare $|\psi(\boldsymbol\theta)\rangle$, apply
$U_b^\dagger$, and measure all qubits in the computational basis.
The all-zero outcome has probability
\begin{equation}\label{eq:overlap_measurement}
p_0(\boldsymbol\theta)
=|\langle0\cdots0|U_b^\dagger
  |\psi(\boldsymbol\theta)\rangle|^2.
\end{equation}
If $N_0$ of $S$ shots return all zeros, set
$\widehat p_0=N_0/S$.

\paragraph{Operator expectation}
Decompose the known scaled matrix into Pauli strings,
\begin{equation}\label{eq:pauli_decomposition}
\overline A=\sum_{\nu=1}^{N_{\rm P}}a_\nu Q_\nu,
\end{equation}
where $N_{\rm P}$ is the number of retained terms,
$a_\nu\in\mathbb R$, and $Q_\nu$ is a tensor product of
single-qubit matrices from $\{I,X,Y,Z\}$. Each measured Pauli
expectation gives a signed average of $+1$ and $-1$ outcomes.
Compatible strings are measured in the same basis setting, and
their estimates are combined as
\[
\widehat D=\sum_{\nu=1}^{N_{\rm P}}
 a_\nu\widehat{\langle Q_\nu\rangle}.
\]

\paragraph{Padding}
Measure the trial state in the computational basis. Outcomes with
indices $N,\ldots,m-1$ correspond to padded amplitudes. If
$N_{\rm pad}$ of $S$ shots have these indices, set
$\widehat p_{\rm pad}=N_{\rm pad}/S$. No padding circuit is needed
when $m=N$.

Because $\overline A$ is SPD and has largest eigenvalue one,
$\lambda_{\min}(\overline A)\leq D\leq1$, where
$\lambda_{\min}$ denotes the smallest eigenvalue. To keep a sampled
denominator positive, define
\[
\widehat D_{\rm safe}
=\min\{1,\max\{\lambda_{\min}(\overline A),\widehat D\}\},
\]
and evaluate
\begin{equation}\label{eq:measured_quantum_objective}
\widehat{\mathcal J}
=-\widehat p_0/\widehat D_{\rm safe}
 +\eta_{\rm pad}\widehat p_{\rm pad}.
\end{equation}
This safeguard and the ratio can bias finite-shot estimates; gate
and readout noise also affect them.}

\subsection{Optimization and classical solution recovery}

{\color{black}A classical simultaneous perturbation stochastic approximation
(SPSA) optimizer~\cite{Spall1992} updates
$\boldsymbol\theta$ from measured objective values at two
perturbed parameter vectors per iteration. The selected final
parameters are denoted by $\boldsymbol\theta^\star$.

{\color{black}We measure $|\psi(\boldsymbol\theta^\star)\rangle$ in all $3^n$
combinations of the single-qubit $X$, $Y$, and $Z$ bases.
Computational-basis probabilities alone do not determine the
relative signs of real amplitudes. State tomography combines the
measurements into an estimated density matrix
$\widehat\rho\in\mathbb C^{m\times m}$~\cite{James2001}.
We take its real symmetric part, clip negative eigenvalues to zero,
and renormalize their sum to one. A normalized eigenvector for the
largest remaining eigenvalue gives the reconstructed real direction
$\widehat\psi\in\mathbb R^m$. This step is a classical
reconstruction and may be imperfect under sampling and noise.

The normalized state specifies a direction but not the magnitude
of the solution. For an ideal real trial direction $\psi$, the magnitude of the energy-minimizing scale could also be estimated from the measured overlap $p_0$ and operator expectation $D$:
$$
|\alpha|
=\frac{\|\widetilde{\mathbf b}\|_2\sqrt{p_0}}{sD},
\qquad \widetilde A_\gamma=s\overline A.
$$
The overlap probability does not determine the sign relative to
a chosen real representative of the state. Since classical
prolongation requires an explicit coarse coefficient vector,
we instead compute the energy-minimizing scale from the
tomographically reconstructed direction $\widehat\psi$ and
remove the padded components:
\begin{align}
\alpha_{\rm rec}
 &=\frac{\widetilde{\mathbf b}^{\top}\widehat\psi}
         {\widehat\psi^\top\widetilde A_\gamma\widehat\psi},
    \label{eq:reconstructed_scaling}\\
\mathbf u_q
 &=R_N(\alpha_{\rm rec}\widehat\psi),
    \label{eq:physical_restriction}
\end{align}
where $R_N=[I_N\;0]\in\mathbb R^{N\times m}$ selects the
first $N$ entries, $I_N$ is the $N\times N$ identity, and the zero
block has $m-N$ columns. An overall sign change of $\widehat\psi$
also changes the sign of $\alpha_{\rm rec}$, leaving
$\mathbf u_q$ unchanged. Both operations are performed on the
classical reconstructed vector, not on the quantum state.

For the Darcy problem, $\mathbf p_H^q=\mathbf u_q$ is the
approximate coarse coefficient vector. Its prolongation gives the
free-node initial pressure for the classical fine-grid correction,
\begin{equation}\label{eq:quantum_fine_initial}
\mathbf p_h^{(0)}=\mathbf g_h+P\mathbf p_H^q.
\end{equation}
The $3^n$ tomography settings make this reconstruction most
appropriate for small coarse registers.}

\section{Quantum Coarse Solve, Lifting, and Fine-Scale Correction}
\label{sec:fine_scale_corr}

This section describes how the quantum coarse solution from
Section~\ref{sec:quantum_solver} enters the classical fine-scale
correction. The quantum computation approximates the
\textcolor{black}{Galerkin coarse solution}; prolongation and
boundary lifting give a fine-grid initial approximation, which a
classical iterative solver then corrects. Thus, the quantum
component serves as a coarse-level solver in the coarse-to-fine
workflow.

In subsurface flow, strong heterogeneity can require fine meshes to
resolve local pressure features. Near-term quantum devices are more
suited to small coarse systems because of qubit, circuit-depth, and
noise limits~\cite{Preskill2018}. This scale separation motivates
the coarse quantum solve followed by classical AMG-preconditioned
fine-grid correction.

The construction follows classical multilevel ideas, including
coarse-grid correction, deflation, and two-level preconditioners
\cite{Nicolaides1987,Saad2003,Smith1996,ToselliWidlund2005,Vuik2006}.
A coarse representation can capture error components that fine-grid
iterations remove slowly. Here, a variational quantum solver supplies
the coarse solution used as an initial approximation for the
classical correction.
\textcolor{black}{
The quantum coarse solve is performed once per linear system to construct
its initial approximation; it does not modify the fine-grid operator or
the AMG preconditioner.
}

Recall that after strongly imposing Dirichlet boundary conditions
and eliminating the prescribed degrees of freedom, the fine-scale
\textcolor{black}{system and the Galerkin coarse system} are
\begin{equation}\label{eq:fine_coarse_systems}
A_h\mathbf{p}_h=\mathbf{b}_h,
\qquad
\textcolor{black}{A_H\mathbf{p}_H^\star=\mathbf{b}_H},
\end{equation}
where
\[
A_h\in\mathbb{R}^{N_h\times N_h},
\qquad
A_H\in\mathbb{R}^{N_H\times N_H},
\]
and
\[
\mathbf{p}_h,\mathbf{b}_h\in\mathbb{R}^{N_h},
\qquad
\textcolor{black}{
\mathbf{p}_H^\star,\mathbf{b}_H\in\mathbb{R}^{N_H}},
\qquad
N_H<N_h.
\]
Here \(A_h\) is the fine-scale stiffness matrix assembled with the
fine permeability field, while
{\color{black}
\[
A_H=P^TA_hP,
\qquad
\mathbf{b}_H=P^T(\mathbf{b}_h-A_h\mathbf{g}_h).
\]
The vector \(\mathbf{g}_h\) contains the free-degree-of-freedom
values of the boundary lifting \(G_h\).
The coarse unknown \(\mathbf{p}_H\) represents the pressure component
after subtracting this lifting, and \(\mathbf{p}_H^\star\) denotes
its exact Galerkin solution.
Thus, the coarse operator is obtained from the fine-scale system;
no independently assembled coarse operator or harmonic-averaged
permeability is used in the proposed solve.
}

Using the variational quantum
\textcolor{black}{coarse solver} described in
Section~\ref{sec:quantum_solver}, we obtain an approximate solution
of the coarse system, denoted by
\textcolor{black}{
\(\mathbf{p}_H^q\in\mathbb{R}^{N_H}\)},
satisfying
\begin{equation}\label{eq:quantum_coarse_solution}
\textcolor{black}{
A_H\mathbf{p}_H^q\approx\mathbf{b}_H.
}
\end{equation}
\textcolor{black}{
This classical coefficient vector is recovered from quantum
measurements by the reconstruction and rescaling procedure described
in Section~\ref{sec:quantum_solver}.
It provides the coarse approximation that will be transferred to
the fine scale.
}

To transfer the quantum-computed coarse approximation to the fine
grid, we use the natural inclusion
\textcolor{black}{
\(I_H^h:V_{H,0}\to V_{h,0}\).
}
Because the fine mesh is a conforming refinement of the coarse mesh,
\(I_H^h\) represents the same coarse function in the fine space.
\textcolor{black}{
Let \(p_H^q\in V_{H,0}\) denote the coarse function associated with
the coefficient vector \(\mathbf{p}_H^q\), so that
\(p_H^q=\sum_{j=1}^{N_H}(\mathbf p_H^q)_j\Phi_j\).
Adding this function to the boundary lifting gives
}
\begin{equation}\label{eq:lifted_function}
\textcolor{black}{
p_h^{(0)}:=G_h+I_H^h p_H^q\in V_{h,g}.
}
\end{equation}
In coefficient form, this inclusion is represented by the
prolongation matrix
\begin{equation}\label{eq:prolongation_matrix}
P\in\mathbb{R}^{N_h\times N_H},
\end{equation}
so that
\begin{equation}\label{eq:lifted_vector}
\textcolor{black}{
\mathbf{p}_h^{(0)}
=
\mathbf{g}_h+P\mathbf{p}_H^q.
}
\end{equation}
Here \(P\) denotes the prolongation between the reduced coarse and
fine coefficient spaces after the same Dirichlet boundary treatment
has been applied on both levels.
\textcolor{black}{
Since \(p_H^q\) vanishes on the Dirichlet boundary, adding its
interpolation to \(G_h\) preserves the prescribed boundary values.
For homogeneous Dirichlet data, we take \(G_h=0\) and
\(\mathbf{g}_h=\mathbf{0}\).
}
The vector \(\mathbf{p}_h^{(0)}\in\mathbb{R}^{N_h}\) contains only
the free-node values of \(p_h^{(0)}\) and is used as the initial
guess for the fine-scale iterative solver.

We perform the fine-scale correction by applying preconditioned
conjugate gradients (PCG) to
\begin{equation}\label{eq:fine_correction_system}
A_h\mathbf{p}_h=\mathbf{b}_h,
\end{equation}
using \(\mathbf{p}_h^{(0)}\) as the initial guess.
The reduced matrix \(A_h\) is SPD, as shown in
Section~\ref{sec:model_problem}. 
We assume that the algebraic multigrid (AMG) preconditioner defines a fixed
SPD operator \(M^{-1}\), so that standard PCG theory
applies to the symmetrically preconditioned operator
\(M^{-1/2}A_hM^{-1/2}\).
\textcolor{black}{
A nonsymmetric or variable AMG cycle would require an appropriate
alternative Krylov method.}

{\color{black}
Let
\begin{equation}\label{eq:initial_residual}
\mathbf{r}_0=\mathbf{b}_h-A_h\mathbf{p}_h^{(0)}
\end{equation}
denote the initial residual associated with the lifted quantum
coarse approximation.
Since the fine-grid solution \(\mathbf{p}_h\) satisfies
\(A_h\mathbf{p}_h=\mathbf{b}_h\), we have
\[
\mathbf{r}_0
=
A_h(\mathbf{p}_h-\mathbf{p}_h^{(0)}).
\]
Thus, the lifted coarse approximation determines the initial error
and residual seen by the fine-grid solver; a smaller initial error
does not alone guarantee fewer iterations.}
\textcolor{black}{
Equivalently, the exact fine-grid correction
\(\mathbf{c}_h=\mathbf{p}_h-\mathbf{p}_h^{(0)}\)
satisfies \(A_h\mathbf{c}_h=\mathbf{r}_0\).
}

{\color{black}Using this preconditioner, the PCG iterates generate a sequence
\(\{\mathbf{p}_h^{(k)}\}\) with residuals
\begin{equation}
\mathbf{r}_k=\mathbf{b}_h-A_h\mathbf{p}_h^{(k)}.
\end{equation}
Formally, the iterates satisfy
\begin{equation}
\mathbf{p}_h^{(k)}
\in
\mathbf{p}_h^{(0)}
+
\mathcal{K}_k(M^{-1}A_h,M^{-1}\mathbf{r}_0),
\end{equation}
where
\[
\mathcal{K}_k(M^{-1}A_h,M^{-1}\mathbf{r}_0)
=
\mathrm{span}
\left\{
M^{-1}\mathbf{r}_0,
(M^{-1}A_h)M^{-1}\mathbf{r}_0,
\dots,
(M^{-1}A_h)^{k-1}M^{-1}\mathbf{r}_0
\right\}.
\]
}
{\color{black}Thus, the quantum coarse solve affects the fine-grid iteration
through the initial residual \(\mathbf{r}_0\) and the initial error
\(\mathbf{e}_0:=\mathbf{p}_h-\mathbf{p}_h^{(0)}\).
To make this dependence explicit, we write

\[
\mathbf{p}_h^{(0)}
=
\mathbf{g}_h+P\mathbf{p}_H^q,
\qquad
\mathbf{p}_H^q
=
\mathbf{p}_H^\star+\boldsymbol{\delta}_H,
\]
}
where
\textcolor{black}{
\(\mathbf{p}_H^\star\) is the exact solution of the Galerkin coarse
system and
\(\boldsymbol{\delta}_H=\mathbf{p}_H^q-\mathbf{p}_H^\star\)
denotes the error in the reconstructed quantum coarse solution.
}
Then
{\color{black}
\[
\mathbf{p}_h-\mathbf{p}_h^{(0)}
=
\left(
\mathbf{p}_h-\mathbf{g}_h-P\mathbf{p}_H^\star
\right)
-
P\boldsymbol{\delta}_H.
\]
The Galerkin construction gives~\cite{Briggs2000}
\[
P^TA_h
\left(
\mathbf{p}_h-\mathbf{g}_h-P\mathbf{p}_H^\star
\right)
=
\mathbf{0}.
\]
Consequently, the two error components are orthogonal in the
energy inner product, and
\begin{equation}\label{eq:initial_energy_decomposition}
\|\mathbf{e}_0\|_{A_h}^2
=
\|\mathbf{p}_h-\mathbf{g}_h-P\mathbf{p}_H^\star\|_{A_h}^2
+
\|\boldsymbol{\delta}_H\|_{A_H}^2,
\end{equation}
where \(\|\mathbf{v}\|_A:=(\mathbf{v}^TA\mathbf{v})^{1/2}\).
The first term represents the error of the exact Galerkin coarse
approximation, while the second represents the additional error
introduced by the approximate quantum coarse solve and classical
reconstruction.
The subsequent AMG-preconditioned CG iteration reduces this initial
fine-grid error according to the convergence behavior of the
preconditioned operator.}

In particular, for a fixed symmetric positive definite preconditioner
\(M^{-1}\), the standard PCG estimate gives
\[
\|\mathbf{e}_k\|_{A_h}
\le
2
\left(
\frac{
\sqrt{\kappa(M^{-1/2}A_hM^{-1/2})}-1
}{
\sqrt{\kappa(M^{-1/2}A_hM^{-1/2})}+1
}
\right)^k
\|\mathbf{e}_0\|_{A_h},
\qquad
\mathbf{e}_k=\mathbf{p}_h-\mathbf{p}_h^{(k)},
\]
where \(\kappa(\cdot)\) denotes the spectral condition number.
Thus, the quantum coarse solve influences the size and structure
of the initial error \(\mathbf{e}_0\), while the AMG preconditioner
influences the worst-case convergence factor through the condition
number of the symmetrically preconditioned operator.
In practice, the observed convergence also depends on how the
initial error is distributed among the modes of the preconditioned
system.

For completeness, the AMG-preconditioned conjugate gradient iteration
is given as follows.
Starting with the initial guess \(\mathbf{p}_h^{(0)}\)
\textcolor{black}{
and the residual \(\mathbf{r}_0\) defined in
\eqref{eq:initial_residual}},
define
\begin{equation}
\mathbf{z}_0=M^{-1}\mathbf{r}_0,
\qquad
\mathbf{d}_0=\mathbf{z}_0.
\end{equation}
Then, for \(k=0,1,2,\dots\), compute
\begin{align}
\alpha_k
&=
\frac{\mathbf{r}_k^\top\mathbf{z}_k}
{\mathbf{d}_k^\top A_h\mathbf{d}_k},
\\
\mathbf{p}_h^{(k+1)}
&=
\mathbf{p}_h^{(k)}+\alpha_k\mathbf{d}_k,
\\
\mathbf{r}_{k+1}
&=
\mathbf{r}_k-\alpha_k A_h\mathbf{d}_k,
\\
\mathbf{z}_{k+1}
&=
M^{-1}\mathbf{r}_{k+1},
\\
\beta_k
&=
\frac{\mathbf{r}_{k+1}^\top\mathbf{z}_{k+1}}
{\mathbf{r}_k^\top\mathbf{z}_k},
\\
\mathbf{d}_{k+1}
&=
\mathbf{z}_{k+1}+\beta_k\mathbf{d}_k.
\end{align}
The iteration terminates when the prescribed relative residual
tolerance is reached. Algorithm~\ref{alg:hybrid_solver} gives the
complete procedure, and Figure~\ref{fig:hybrid_workflow} summarizes the complete hybrid
quantum--classical coarse-to-fine workflow.

For the time-dependent problem, the same initialization and
correction are applied at each time step to its SPD time-discrete
system, using that step's operator and right-hand side.

\begin{algorithm}[!htbp]
\caption{Hybrid quantum--classical coarse-to-fine Darcy solver}
\label{alg:hybrid_solver}
\begin{enumerate}

\item \textbf{Input:} Reduced fine-grid system
$A_h\mathbf p_h=\mathbf b_h$ with $\mathbf b_h\ne\mathbf 0$;
prolongation $P\in\mathbb R^{N_h\times N_H}$; boundary lifting
$G_h$ and its free-node vector $\mathbf g_h$; fixed SPD AMG
preconditioner $M^{-1}$; quantum solver settings; PCG tolerance
$\varepsilon_{\mathrm{pcg}}$ and iteration limit $k_{\max}$.

\item \textbf{Form the coarse system:}
\[
A_H=P^TA_hP,\qquad
\mathbf b_H=P^T(\mathbf b_h-A_h\mathbf g_h).
\]
{\color{black}
\item \textbf{Solve and reconstruct on the coarse grid:}
If $\mathbf b_H=\mathbf 0$, set $\mathbf p_H^q=\mathbf 0$.
Otherwise, pad and scale $(A_H,\mathbf b_H)$ as in
\eqref{eq:padded_system}-\eqref{eq:scaled_quantum_operator}.
Prepare the right-hand-side state with $U_b$ and trial states
$|\psi(\boldsymbol\theta)\rangle$. Run SPSA on
\eqref{eq:measured_quantum_objective} and select
$\boldsymbol\theta^\star$. At these parameters, perform tomography
to reconstruct $\widehat\psi$. Compute $\alpha_{\rm rec}$ using
\eqref{eq:reconstructed_scaling} and set
\[
\mathbf p_H^q=R_{N_H}(\alpha_{\rm rec}\widehat\psi),
\]
where $R_{N_H}$ keeps the first $N_H$ components.

\item \textbf{Initialize the fine-grid solve:}
\[
\mathbf p_h^{(0)}=\mathbf g_h+P\mathbf p_H^q,\qquad
\mathbf r_0=\mathbf b_h-A_h\mathbf p_h^{(0)},\qquad
\mathbf z_0=M^{-1}\mathbf r_0,\quad \mathbf d_0=\mathbf z_0.
\]
Set $k=0$.}

\item \textbf{Apply AMG-PCG:} While
$\|\mathbf r_k\|_2/\|\mathbf b_h\|_2>
\varepsilon_{\mathrm{pcg}}$ and $k<k_{\max}$, compute
\[
\begin{aligned}
\alpha_k
 &=\frac{\mathbf r_k^T\mathbf z_k}
         {\mathbf d_k^TA_h\mathbf d_k},\\
\mathbf p_h^{(k+1)}
 &=\mathbf p_h^{(k)}+\alpha_k\mathbf d_k,\\
\mathbf r_{k+1}
 &=\mathbf r_k-\alpha_kA_h\mathbf d_k,\\
\mathbf z_{k+1}
 &=M^{-1}\mathbf r_{k+1},\\
\beta_k
 &=\frac{\mathbf r_{k+1}^T\mathbf z_{k+1}}
         {\mathbf r_k^T\mathbf z_k},\\
\mathbf d_{k+1}
 &=\mathbf z_{k+1}+\beta_k\mathbf d_k.
\end{aligned}
\]
Then set $k\leftarrow k+1$.

\item \textbf{Output:} Return the corrected free-node vector
$\mathbf p_h^{\mathrm{corr}}=\mathbf p_h^{(k)}$ and restore the
Dirichlet values to form $p_h^{\mathrm{corr}}$. Record $k$,
$\|\mathbf r_0\|_2$, and $\|\mathbf r_k\|_2$.

\end{enumerate}
\end{algorithm}

\begin{figure}[!htbp]
\centering
\includegraphics[width=\linewidth]{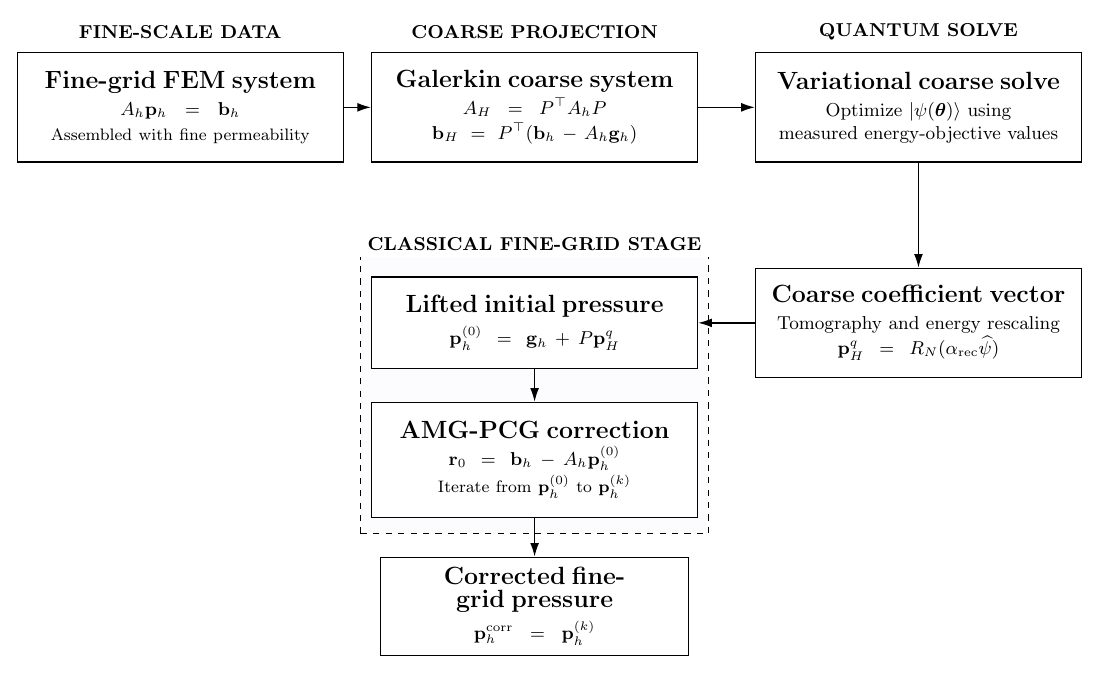}
\caption{Hybrid quantum-classical coarse-to-fine workflow.}
\label{fig:hybrid_workflow}
\end{figure}

\section{Numerical Experiments}
\label{sec:numerical_experiments}

We assess the hybrid coarse-to-fine framework on steady and
time-dependent heterogeneous Darcy problems. We examine whether
the reconstructed coarse solution provides a useful initial
guess for the fine-grid solve. We report coarse-solution
accuracy, initial and corrected pressure errors, AMG-PCG
iteration counts, and the sampling requirements of the
coarse solve.

The reference solution $p_h$ is computed from the same
fine-grid finite element system used for correction.
Bold symbols denote coefficient vectors, while the
corresponding plain symbols denote finite element functions.
Pressure vectors contain only the values at nodes where
pressure is unknown; prescribed Dirichlet values are included
when forming the full pressure functions.

For finite element functions $u$ and a nonzero reference $v$,
we define the relative $L^2(\Omega)$ error by
\begin{equation}
\label{eq:relative_error}
E_{\mathrm{rel}}(u,v)
:=
\frac{\|u-v\|_{L^2(\Omega)}}{\|v\|_{L^2(\Omega)}}.
\end{equation}
We evaluate these norms using the assembled fine-grid mass
matrix and full nodal vectors.

For the steady problem, $\mathbf p_H^q$ contains the
reconstructed coarse coefficients before adding the boundary
lift. Its finite element function $p_H^q$ has zero Dirichlet
values. Adding the boundary lift $G_H$ gives the reconstructed
coarse solution
\[
\widetilde p_H^q=G_H+p_H^q.
\]
Let $\mathbf p_H^\star$ be the exact solution of
$A_H\mathbf p_H^\star=\mathbf b_H$. We assess the reconstructed
coefficients using the relative coarse residual and
relative coefficient error:
\begin{equation}
\label{eq:ex1_coarse_metrics}
\rho_H
:=
\frac{\|A_H\mathbf p_H^q-\mathbf b_H\|_2}
     {\|\mathbf b_H\|_2},
\qquad
E_{H,2}
:=
\frac{\|\mathbf p_H^q-\mathbf p_H^\star\|_2}
     {\|\mathbf p_H^\star\|_2}.
\end{equation}
The residual $\rho_H$ measures how well the reconstructed
vector satisfies the coarse system. The coefficient error
$E_{H,2}$ measures its relative Euclidean distance from the
exact Galerkin vector.

For the steady fine-grid solve, we compare three initializations:
\begin{equation}
\label{eq:fine_initializations}
\mathbf p_h^{(0)}=
\begin{cases}
\mathbf 0,
& \text{zero initialization},\\
\mathbf g_h+P\mathbf p_H^\star,
& \text{exact Galerkin initialization},\\
\mathbf g_h+P\mathbf p_H^q,
& \text{reconstructed coarse initialization}.
\end{cases}
\end{equation}
Here, $P$ prolongs the coarse coefficients to the fine grid,
and $\mathbf g_h$ contains the values of $I_H^hG_H$ at the
nodes where pressure is unknown. For reconstructed coarse
initialization, the corresponding fine-grid function is
$p_h^{(0)}=I_H^h\widetilde p_H^q$. All three initializations
retain the prescribed Dirichlet values.

For an unknown-pressure vector $\mathbf v\in\mathbb R^{N_h}$,
we define the relative fine-grid residual and energy error by
\begin{equation}
\label{eq:ex1_fine_metrics}
\rho_h(\mathbf v)
:=
\frac{\|\mathbf b_h-A_h\mathbf v\|_2}
     {\|\mathbf b_h\|_2},
\qquad
E_A(\mathbf v)
:=
\left(
\frac{(\mathbf v-\mathbf p_h)^T A_h(\mathbf v-\mathbf p_h)}
     {\mathbf p_h^T A_h\mathbf p_h}
\right)^{1/2},
\end{equation}
where $\mathbf p_h$ contains the unknown-node coefficients
of the reference solution.

For each initialization in \eqref{eq:fine_initializations},
we report
\[
\rho_0:=\rho_h(\mathbf p_h^{(0)}),
\qquad
E_{\mathrm{rel}}^{(0)}
:=E_{\mathrm{rel}}(p_h^{(0)},p_h),
\qquad
E_A^{(0)}:=E_A(\mathbf p_h^{(0)}).
\]
After AMG-PCG correction, we report the iteration count and
\[
\rho_{\mathrm{corr}}
:=\rho_h(\mathbf p_h^{\mathrm{corr}}),
\qquad
E_{\mathrm{rel}}^{\mathrm{corr}}
:=E_{\mathrm{rel}}(p_h^{\mathrm{corr}},p_h).
\]
In the time-dependent experiment, residuals and $L^2$ errors
are evaluated at each time step using the backward Euler
system and a reference computed on the same spatial mesh
and time grid. All reported relative quantities have nonzero
denominators.

{\color{black}
Both experiments are implemented in Python. We use NumPy
and SciPy for finite element assembly and sparse linear
algebra. Fine-grid correction uses SciPy's conjugate-gradient
solver with a PyAMG smoothed-aggregation preconditioner
(AMG-PCG). We use the default V-cycle settings and construct
the AMG hierarchy once for each fixed system matrix.

The quantum circuits are implemented with Qiskit and simulated
with Qiskit Aer~\cite{qiskit2021}. Objective evaluations and
tomography use finite-shot measurements, with and without
modeled gate and readout noise. Ideal statevector calculations
provide reference results and the fixed circuit parameters
used in the shot-count study. All quantum computations
are simulated.

\paragraph{Measurement conditions}
\label{par:measurement_conditions}
Each measurement setting uses a finite number of shots, denoted
by $S$. One shot is one circuit execution followed by measurement.
Because the number of shots is limited, the measured results
can vary even with ideal gates and readout. We also consider
synthetic noise and device-modeled noise.

For synthetic gate noise, the simulator changes the state using
the depolarizing channel
\begin{equation}
\mathcal D_{\lambda}^{(n_{\mathrm{gate}})}(\rho)
=(1-\lambda)\rho
+\lambda\frac{I_{2^{n_{\mathrm{gate}}}}}{2^{n_{\mathrm{gate}}}}.
\end{equation}
Here, $\rho$ describes the state of the qubits affected by
the gate. We set $n_{\mathrm{gate}}=1$ for a one-qubit gate
and $n_{\mathrm{gate}}=2$ for a two-qubit gate. The matrix $I_{2^{n_{\mathrm{gate}}}}$ is the
$2^{n_{\mathrm{gate}}}\times 2^{n_{\mathrm{gate}}}$ identity
matrix
, and
$I_{2^{n_{\mathrm{gate}}}}/2^{n_{\mathrm{gate}}}$ represents
a fully random state. The parameter $\lambda$ sets how much
of this random state is mixed into $\rho$
\cite{qiskitAerDepolarizing}.
We apply this noise after the compiled single-qubit
Pauli-$X$ and square-root-of-$X$ gates with strength
$\lambda_1$, and after the two-qubit CNOT gates with
strength $\lambda_2$. Here, $X$ denotes the quantum
bit-flip operation.

For synthetic readout noise, each measured bit is independently
recorded as the opposite value with probability
$\epsilon_{\mathrm{ro}}$:
\begin{equation}
\Pr(\text{record }1\mid 0)
=\Pr(\text{record }0\mid 1)
=\epsilon_{\mathrm{ro}}.
\end{equation}
Thus, a $0$ may be recorded as $1$, or a $1$ as $0$.
This changes the recorded result rather than the state
before measurement~\cite{qiskitAerReadout}.
The mild and moderate synthetic noise settings use
different values of $(\lambda_1,\lambda_2,\epsilon_{\mathrm{ro}})$,
specified in the experiments below.

For device-modeled noise, Qiskit Aer uses gate and readout
errors from the stored \texttt{FakeManilaV2} calibration
snapshot~\cite{qiskitFakeProvider}. These errors come from
the stored device calibration data rather than the
synthetic noise settings above.}

\subsection{Experiment 1: Two-Dimensional Heterogeneous Darcy Flow}
\label{sec:ex1}

In this experiment, we test whether the reconstructed coarse
solution improves the initial fine-grid solution
$\mathbf p_h^{(0)}$ and reduces AMG-PCG iterations.
We also examine how shot count and modeled noise
affect coarse-solution accuracy.

\paragraph{Problem setup}
We consider the steady Darcy problem~\eqref{eq:model} on
$\Omega=(0,1)^2$. we impose Dirichlet boundary conditions on the left and the right boundaries ($p(0,y)=1$ and $p(1,y)=0$), together with homogeneous no-flow
conditions on the top and bottom boundaries and set the source
term to $f(x,y)=10\sin(\pi x)\sin(\pi y)$. The fine and Galerkin coarse meshes contain 8,192 and 32 elements, respectively. Using linear finite elements and eliminating the prescribed Dirichlet boundary values give $N_h=4{,}095$ and $N_H=15$ degrees of
freedom, respectively. We use the
heterogeneous fine-scale permeability $K_h=\kappa_h I$ shown in
Figure~\ref{fig:fine_reference_ex1}. Its fine-scale
variations make it challenging for the coarse solution to
capture the local behavior of the pressure field.
{\color{black}Following \eqref{eq:padded_system}, we augment the 15 coarse
degrees of freedom with one auxiliary unknown, giving a
16-dimensional padded system represented on four qubits.
We set $\gamma=\lambda_{\max}(A_H)$, use a six-layer ansatz
with 24 parameters, and choose the padding-penalty weight
$\eta_{\mathrm{pad}}=1$.}

We present the results in two parts. First, we assess the accuracy of the reconstructed coarse solutions under finite-shot sampling and modeled noise, then test their effect on the initial fine-grid solution and classical solver (AMG-PCG) iteration count. Second, we study how the number of measurement shots affects coarse-solution accuracy while holding the circuit parameters fixed. 

{\color{black}\paragraph{i) Coarse-solution accuracy and its effect on the fine-grid solve} In the first study, we test whether a reconstructed coarse solution from simulated
finite-shot measurements provides a useful initial solution, $p_h^{(0)}$, for
the fine-grid solve. We consider three measurement conditions defined above
(p.~\pageref{par:measurement_conditions}): finite-shot
sampling without gate or readout noise, mild synthetic
noise, and device-modeled noise. For mild synthetic noise,
we use $(\lambda_1,\lambda_2,\epsilon_{\mathrm{ro}})
=(10^{-4},10^{-3},0.005)$.

For each condition, we perform five independent SPSA runs with
200 updates and 2,048 shots per objective circuit. In each run,
we select a candidate using fresh measurements with 4,096 shots
per circuit and reconstruct its state using all 81 product-Pauli
bases with 4,096 shots per basis.

Table~\ref{tab:ex1_quantum_accuracy} reports the relative
coarse residual $\rho_H$ and coefficient error $E_{H,2}$
defined in \eqref{eq:ex1_coarse_metrics}.
Device-modeled noise gives the largest mean values of both
measures. Finite-shot sampling and mild noise give similar
coefficient errors, while mild noise gives a lower mean
coarse residual.}

\begin{table}[!htbp]
\centering
\begin{tabular}{lcc}
\hline
Measurement condition & $\rho_H$ & $E_{H,2}$ \\
\hline
Finite shots & $0.721\pm0.115$ & $0.267\pm0.061$ \\
Mild noise & $0.590\pm0.100$ & $0.282\pm0.030$ \\
Device-modeled noise & $0.797\pm0.103$ & $0.410\pm0.048$ \\
\hline
\end{tabular}
\caption{Experiment 1. {\color{black}Accuracy of reconstructed coarse solutions.
Entries are means and sample standard deviations over five independent
SPSA runs per measurement condition.}}
\label{tab:ex1_quantum_accuracy}
\end{table}

{\color{black}Table~\ref{tab:ex1_fine_correction} compares the initializations
defined in \eqref{eq:fine_initializations}. We report the
initial residual $\rho_0$, relative $L^2$ error
$E_{\mathrm{rel}}^{(0)}$, and energy error $E_A^{(0)}$,
using the definitions in \eqref{eq:relative_error}
and \eqref{eq:ex1_fine_metrics}. Under all three measurement conditions, reconstructed coarse
initialization gives smaller $E_{\mathrm{rel}}^{(0)}$ and
$E_A^{(0)}$ than zero initialization. Exact Galerkin
initialization gives the smallest errors. Both types of
coarse initialization give larger $\rho_0$ than zero
initialization. After correction, all cases have comparable
$E_{\mathrm{rel}}^{\mathrm{corr}}$ and require 26 AMG-PCG
iterations. Thus, in this experiment, we observe smaller
initial pressure errors without a reduction in AMG-PCG
iterations.}

\begin{table}[!htbp]
\centering
\small
\setlength{\tabcolsep}{2pt}
\resizebox{\linewidth}{!}{%
\begin{tabular}{@{}lccccc@{}}
\hline
Initialization & $\rho_0$ & $E_{\mathrm{rel}}^{(0)}$
& $E_A^{(0)}$ & \shortstack{AMG-PCG\\iterations}
& \shortstack{$E_{\mathrm{rel}}^{\mathrm{corr}}$\\$(\times10^{-13})$} \\
\hline
Zero & $1.000$ & $1.000$ & $1.000$ & 26 & $4.42$ \\
Exact Galerkin
& $28.423$ & $0.573$ & $0.758$ & 26 & $3.04$ \\
Finite shots
& $26.184\pm0.785$ & $0.650\pm0.022$ & $0.800\pm0.012$
& 26 & $3.29\pm0.26$ \\
Mild noise
& $26.499\pm0.716$ & $0.649\pm0.016$ & $0.795\pm0.008$
& 26 & $3.26\pm0.23$ \\
Device-modeled noise
& $25.044\pm1.084$ & $0.704\pm0.023$ & $0.823\pm0.012$
& 26 & $3.55\pm0.23$ \\
\hline
\end{tabular}
}
\caption{Experiment 1. {\color{black}Initial residual $\rho_0$, relative
$L^2$ error $E_{\mathrm{rel}}^{(0)}$, energy error $E_A^{(0)}$,
AMG-PCG iterations, and corrected relative $L^2$ error
$E_{\mathrm{rel}}^{\mathrm{corr}}$. The three measurement-condition
rows use reconstructed coarse initialization. Entries with
$\pm$ show means and sample standard deviations over five
independent runs.} }
\label{tab:ex1_fine_correction}
\end{table}

To illustrate the pressure fields, we select the device-modeled noise run
with the median coarse residual. Figure~\ref{fig:pressure_comparison} displays $\kappa_h$, $p_h$,
$\widetilde p_H^q$, $p_h^{(0)}$, $p_h^{\mathrm{corr}}$, and the
pressures along $y=0.5$. The lifted quantum solution
captures the large-scale pressure structure, and classical correction
recovers the fine-grid reference solution to solver accuracy.

\begin{figure}[!htbp]
\centering
\begin{subfigure}[t]{0.32\linewidth}
  \centering
  \includegraphics[width=\linewidth]{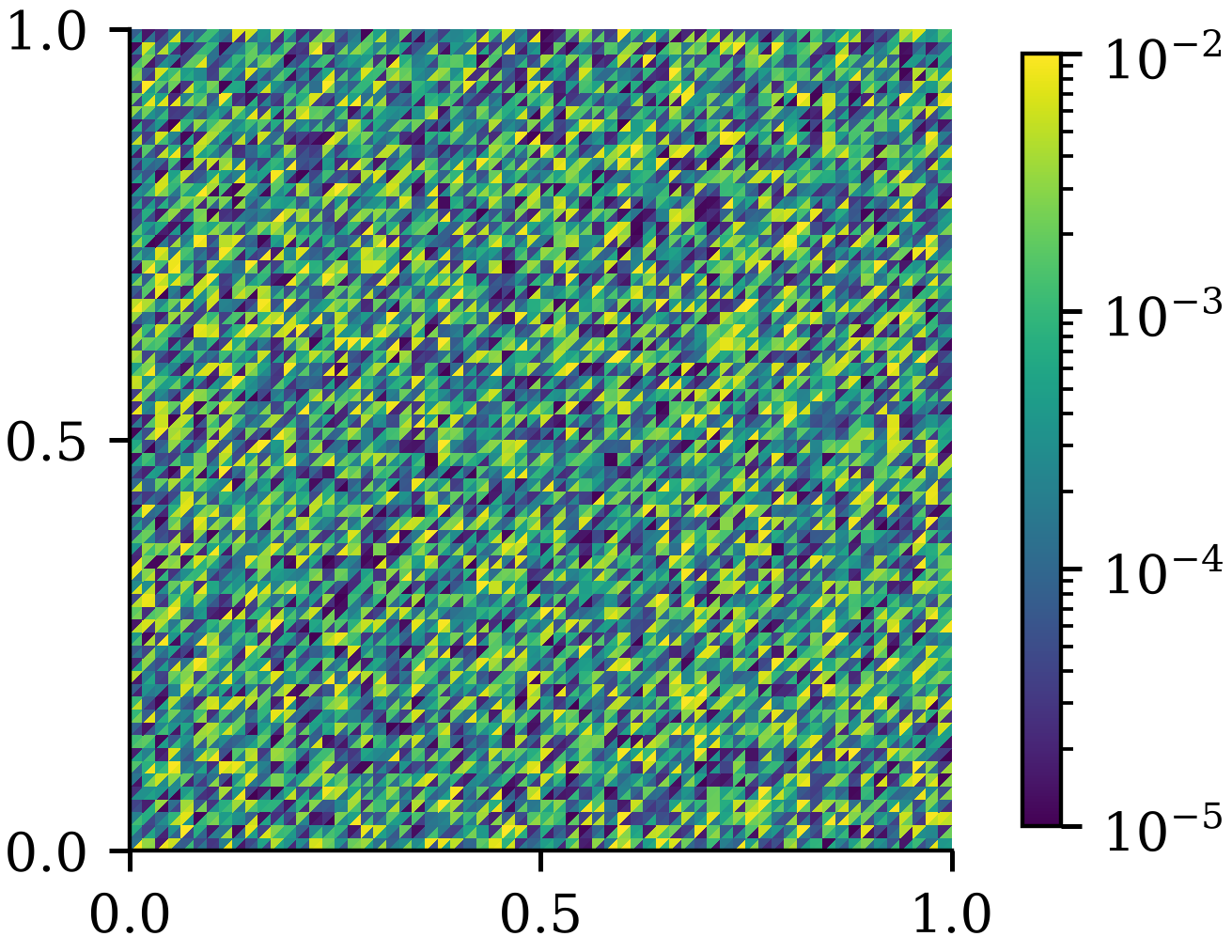}
  \caption{{\color{black}Permeability $\kappa_h$.}}
  \label{fig:fine_reference_ex1}
\end{subfigure}
\hfill
\begin{subfigure}[t]{0.32\linewidth}
  \centering
  \includegraphics[width=\linewidth]{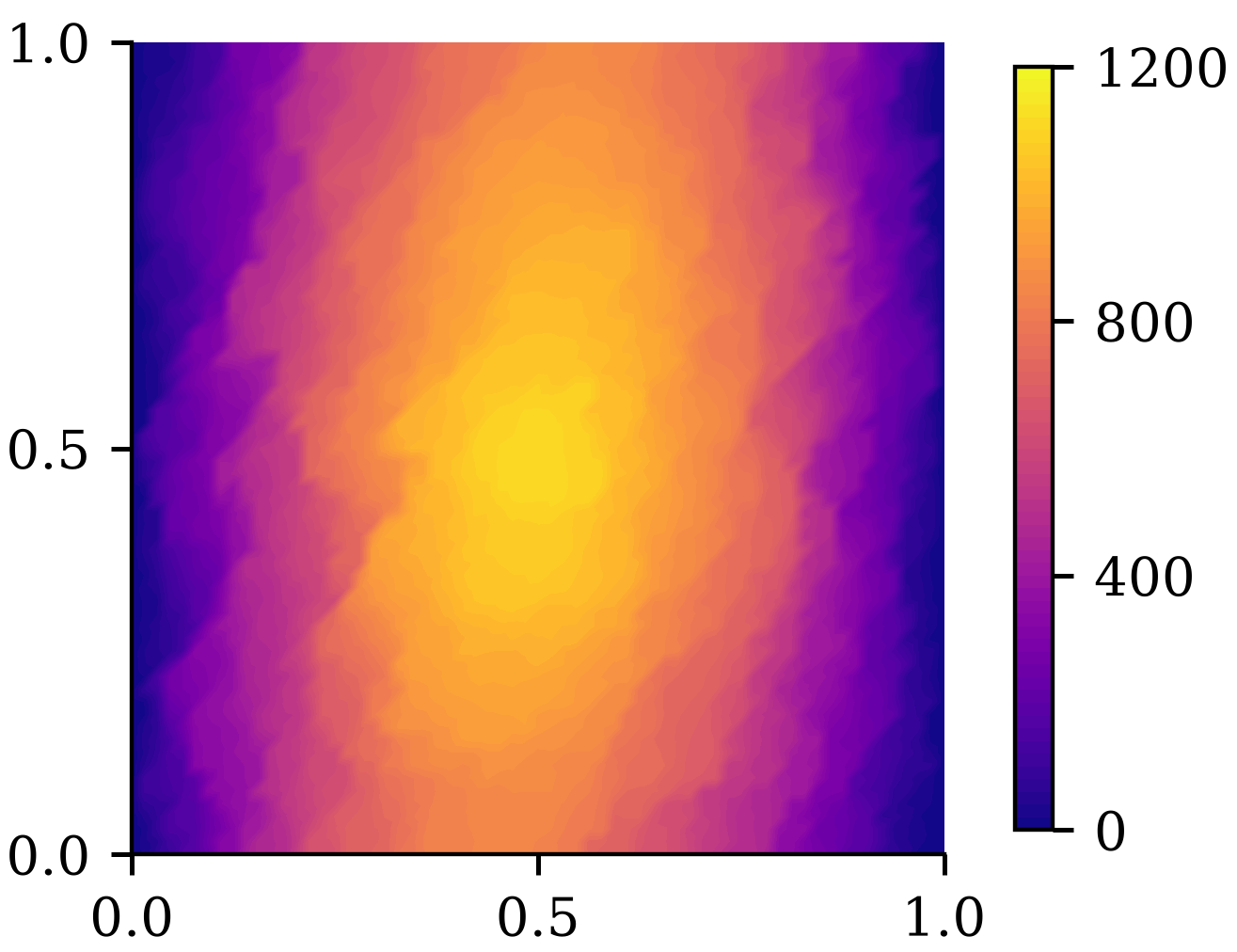}
  \caption{Reference pressure $p_h$.}
  \label{fig:fine_p_ex1}
\end{subfigure}
\hfill
\begin{subfigure}[t]{0.32\linewidth}
  \centering
  \includegraphics[width=\linewidth]{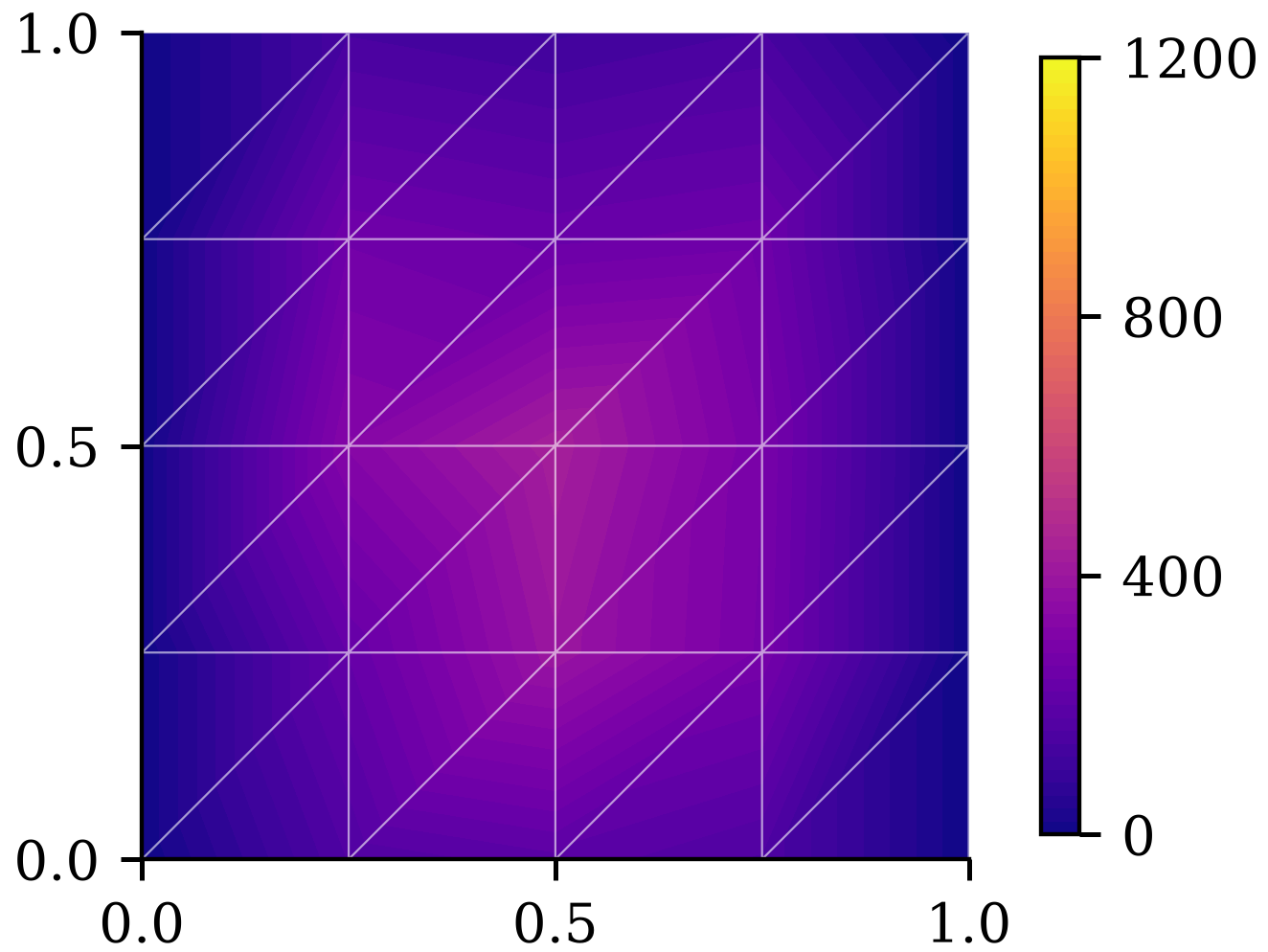}
  \caption{Coarse pressure $\widetilde p_H^q$.}
  \label{fig:coarse_reference_p_ex1}
\end{subfigure}

\par\medskip
\begin{subfigure}[t]{0.32\linewidth}
  \centering
  \includegraphics[width=\linewidth]{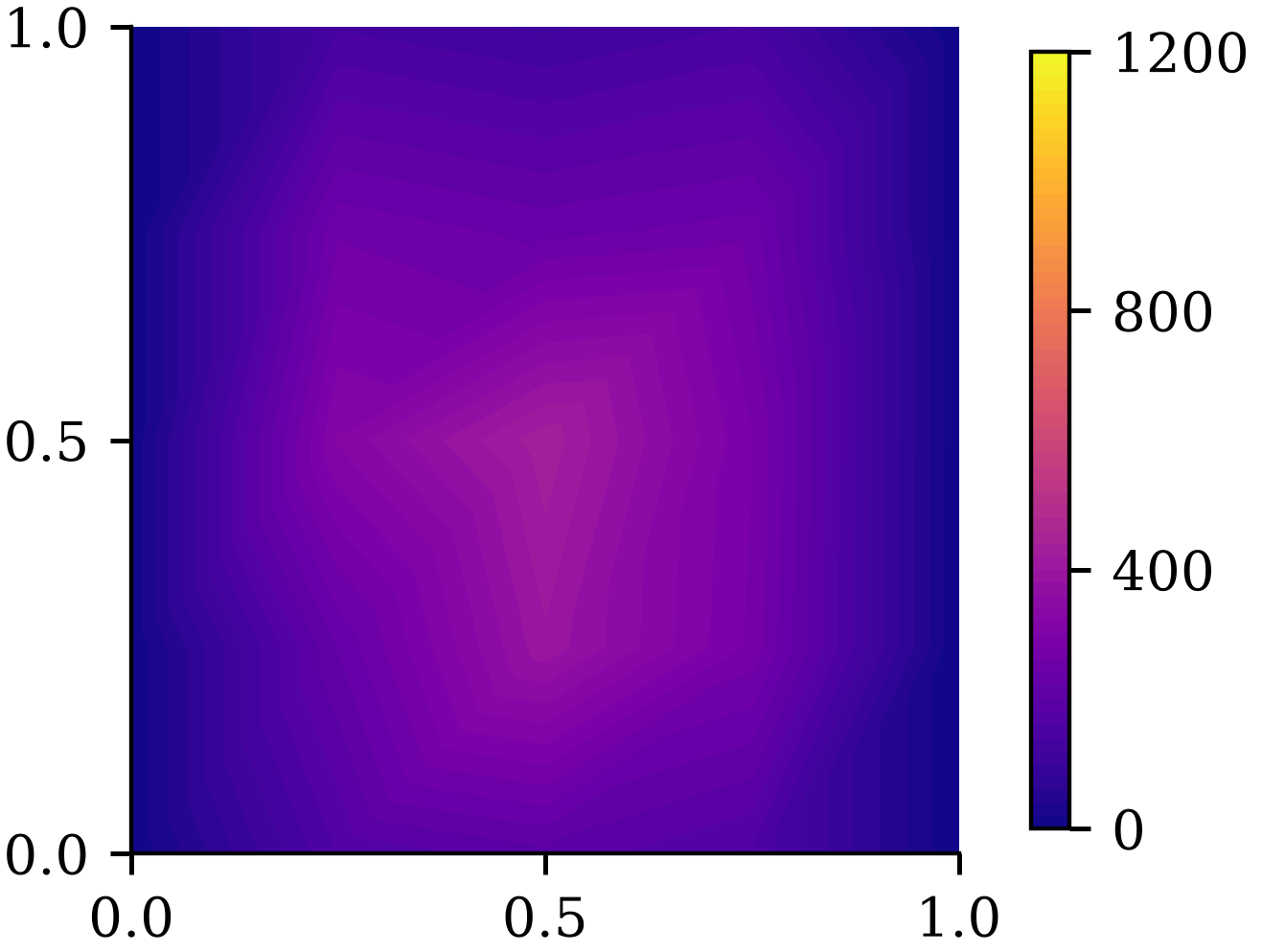}
  \caption{Lifted initial pressure $p_h^{(0)}$.}
  \label{fig:ex1_lifted}
\end{subfigure}
\hfill
\begin{subfigure}[t]{0.32\linewidth}
  \centering
  \includegraphics[width=\linewidth]{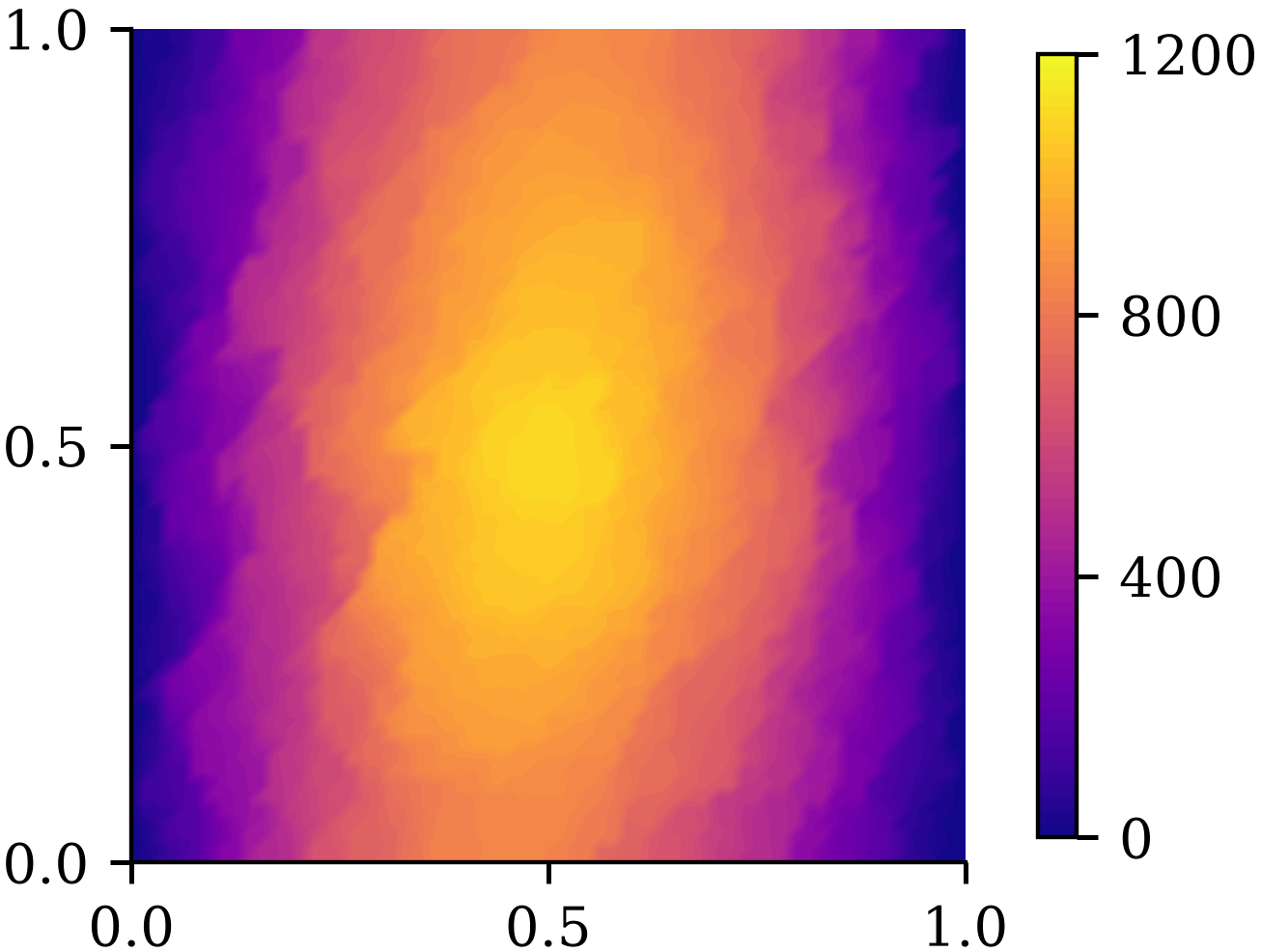}
  \caption{Corrected pressure $p_h^{\mathrm{corr}}$.}
  \label{fig:fine_corrected_p_ex1}
\end{subfigure}
\hfill
\begin{subfigure}[t]{0.32\linewidth}
  \centering
  \includegraphics[width=\linewidth]{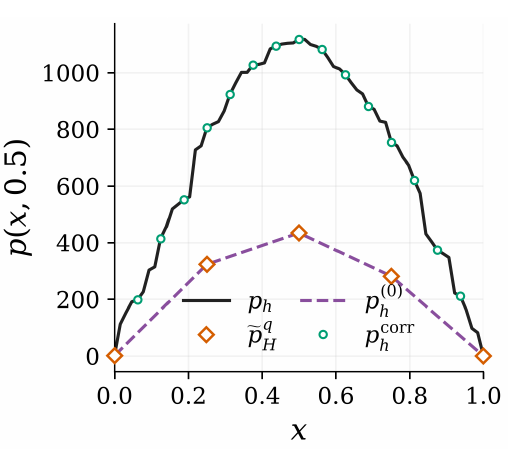}
  \caption{Pressure slice at $y=0.5$.}
  \label{fig:fine_corrected_ex1}
\end{subfigure}
\caption{Experiment 1. fields for the device-modeled noise run with the median
coarse residual: (a) $\kappa_h$, (b) $p_h$,
(c) $\widetilde p_H^q$, (d) $p_h^{(0)}=I_H^h\widetilde p_H^q$,
and (e) $p_h^{\mathrm{corr}}$. Panel (f) compares these pressures
along $y=0.5$: $p_h$ (solid line), $\widetilde p_H^q$ (diamonds at
coarse nodes), $p_h^{(0)}$ (dashed line), and
$p_h^{\mathrm{corr}}$ (circles).}
\label{fig:pressure_comparison}
\end{figure}

{\color{black}
Table~\ref{tab:ex1_sampling_cost} reports measurement costs
for objective evaluation, tomography, and a complete SPSA run.
Each objective uses 22 Pauli groups, an overlap circuit,
and a padding circuit. The compiled standalone $U_b$ has
depth 54 and 23 CNOT gates; $U_b^\dagger$ is executed on
every overlap-measurement shot. The complete run includes
gain calibration, 200 SPSA updates, checkpoint measurements,
candidate selection, and final readout.
}
\begin{table}[!htbp]
\centering
\small
\begin{tabular}{lrrr}
\hline
Measurement & Circuit executions & Shots per circuit & Total shots \\
\hline
One SPSA objective evaluation & 24 & 2,048 & 49,152 \\
Final tomography & 81 & 4,096 & 331,776 \\
\hline
Complete measured SPSA run & 10,521 & - & 22,302,720 \\
\hline
\end{tabular}
\caption{Experiment 1. {\color{black}Simulated circuit and shot counts for the
measured SPSA runs in the first study. Shots per circuit vary within
the complete run.}}
\label{tab:ex1_sampling_cost}
\end{table}

\paragraph{ii) Effect of shot count}

The first study used fixed shot counts for each measurement
stage. We now vary the shot count to examine how it affects
the accuracy of the reconstructed coarse solution.
To separate measurement effects from differences between
SPSA runs, we fix the circuit parameters at values obtained
from ideal statevector optimization.

For each shot count and measurement condition, we repeat
the objective evaluation and state reconstruction ten times
using fresh measurement samples. In each repetition, we use
the selected shot count for every objective circuit and
tomography basis. We use the noise models as previously described
(p.~\pageref{par:measurement_conditions}).
Alongside the three conditions used in the first study,
we include moderate synthetic noise with
$(\lambda_1,\lambda_2,\epsilon_{\mathrm{ro}})
=(10^{-3},10^{-2},0.02)$.

Figure~\ref{fig:ex1_shot_study} shows that the mean coarse
residual decreases with increasing shot count under all four
measurement conditions. Mild noise gives results close to
finite-shot sampling, while moderate synthetic noise and
device-modeled noise give larger residuals. Device-modeled
noise gives the largest mean residual at every shot count.
Thus, increasing the shot count reduces the coarse residual,
but differences between measurement conditions persist
over the range tested.

\begin{figure}[!h]
\centering
\includegraphics[width=0.6\linewidth]{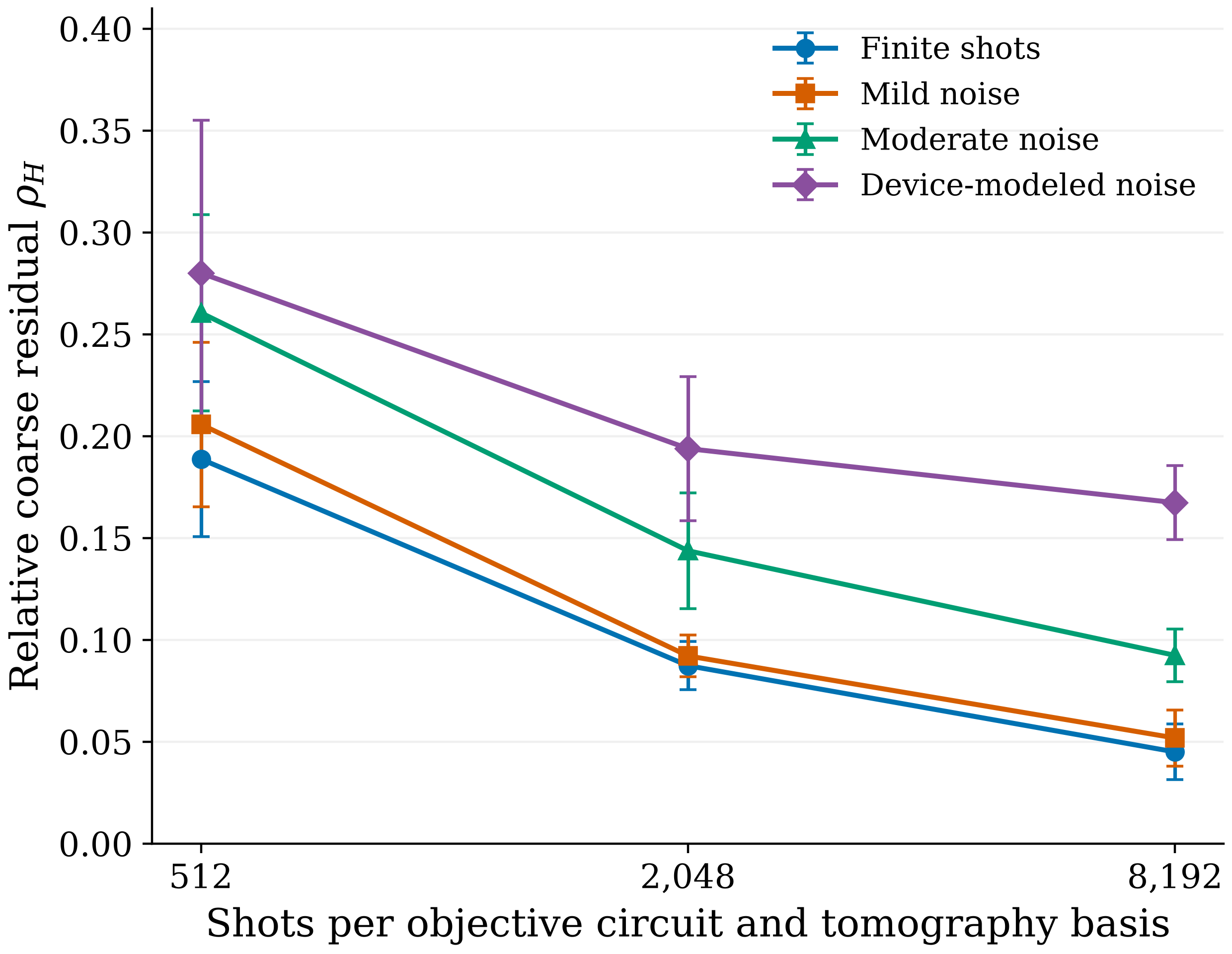}
\caption{{\color{black}Experiment 1. Relative coarse residual $\rho_H$ versus shots per
objective circuit and tomography basis. Circuit parameters
are fixed at values obtained from ideal statevector
optimization. Points show means and error bars show sample
standard deviations over ten independent measurement
repetitions per setting.}}
\label{fig:ex1_shot_study}
\end{figure}

{\color{black}
\subsection{Experiment 2: Time-Dependent Darcy Flow with SPE10 Permeability}
\label{sec:ex2}

In this experiment, we apply the coarse-to-fine correction framework
to a time-dependent Darcy problem with a heterogeneous SPE10-based
permeability field~\cite{christie2001tenth}. Time discretization produces a linear system at each time step.
We approximately solve the corresponding coarse system using
the quantum circuit and reconstruct the coarse solution from
measurements. We then prolongate this solution to the fine grid
and use it to initialize the classical AMG-PCG solve.

We assess the accuracy
of the reconstructed coarse solution, the fidelity of its associated quantum state, and the pressure error after fine-grid correction. We also examine
whether reusing circuit parameters reduces sampling requirements
and compare AMG-PCG iteration counts
for initialization from the reconstructed coarse solution
and the exact Galerkin coarse solution.

\paragraph{Problem setup}
We consider the pressure field
$p=p(\mathbf x,t):\Omega\times[0,T]\to\mathbb R$, which satisfies
\begin{equation}
\frac{\partial p}{\partial t}
-\nabla\cdot\bigl(K_h(\mathbf x)\nabla p\bigr)
=f(\mathbf x,t)
\qquad\text{in }\Omega\times(0,T].
\label{eq:ex2_model}
\end{equation}
Here we set $\Omega=(0,L_x)\times(0,60)$, $L_x=220$, and $T=1$, and  we impose Dirichlet boundary conditions on the left and the right boundaries 
($p(0,y,t)=1$ and $p(L_x,y,t)=0$), together with homogeneous
no-flow conditions on the upper and lower boundaries. The initial condition is
$p(\mathbf x,0)=1-x/L_x$. We set the source term to
\begin{equation}
\begin{split}
f(\mathbf x,t)&=a(t)\bigl(G_{0.3}(\mathbf x)-G_{0.7}(\mathbf x)\bigr),\\
G_c(\mathbf x)&=\exp\left(-\frac{(x/L_x-c)^2+(y/60-0.5)^2}{0.01}\right),
\qquad
 a(t)=\begin{cases}10,&0<t\leq0.5,\\20,&0.5<t\leq1.\end{cases}
\end{split}
\label{eq:ex2_forcing}
\end{equation}
The two Gaussian terms represent a localized source-sink pair. Their amplitudes
double after $t=0.5$, producing a transient pressure response with fixed
permeability and boundary conditions. 

We use the heterogeneous permeability
$K_h(\mathbf x)=\kappa_h(\mathbf x)I$ obtained from the top
two-dimensional layer of the SPE10 dataset~\cite{christie2001tenth}
and shown in Figure~\ref{fig:ex2_permeability}.
The fine and Galerkin coarse meshes contain 26,400 and
32 elements, respectively. Using linear finite elements
and eliminating the prescribed Dirichlet boundary values
give $N_h=13{,}359$ and $N_H=15$ degrees of freedom,
respectively. {\color{black}Following \eqref{eq:padded_system}, we augment the 15 coarse
degrees of freedom with one auxiliary unknown, giving a
16-dimensional padded system represented on four qubits.
As in Experiment~\ref{sec:ex1}, we use a six-layer ansatz
with 24 parameters.}
\begin{figure}[!h]
\centering
\includegraphics[width=.85\linewidth]{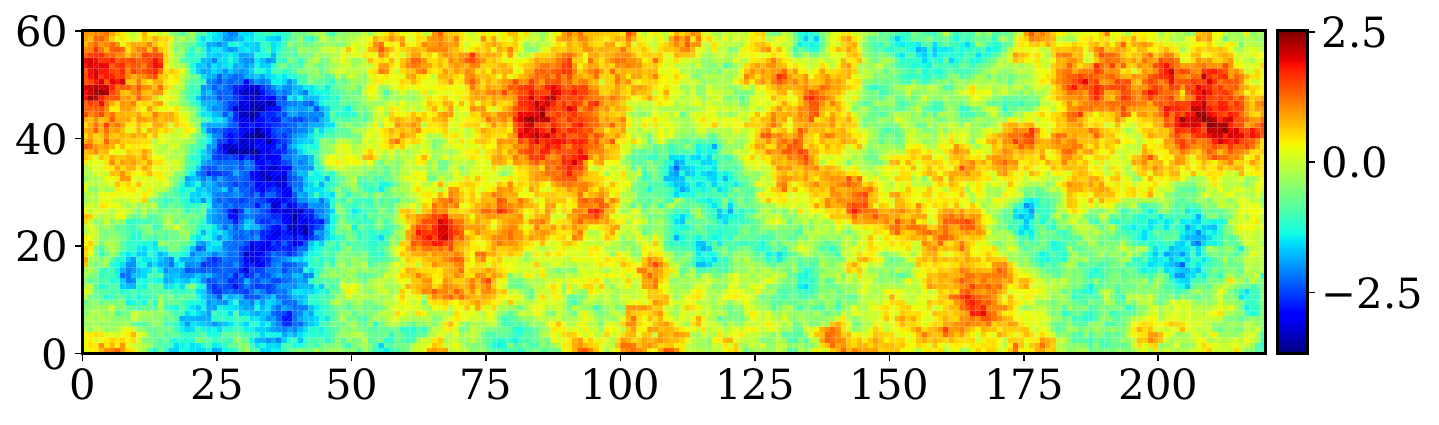}
\caption{{\color{black}Experiment 2. fine permeability
$\kappa_h(\mathbf{x})$
from a two-dimensional top slice of the SPE10 dataset; colors show
$\log_{10}\kappa_h$.}}
\label{fig:ex2_permeability}
\end{figure}

\paragraph{Time discretization and coarse solution}
We use backward Euler time discretization with $\Delta t=0.05$,
giving $N_t=T/\Delta t=20$ time steps and $t_i=i\Delta t$.
Let $M_h^{\mathrm{full}}$ and $A_h^{\mathrm{full}}$ denote the
fine-grid mass and Darcy stiffness matrices before eliminating
the prescribed pressure values. Let $F$ index the nodes where
pressure is unknown and $D$ index the nodes with prescribed
Dirichlet values. At time $t_i$, the fine-grid system is
\begin{equation}
\begin{split}
B_h\mathbf p_h^i
&=\mathbf b_h^i,
\qquad
B_h=
(M_h^{\mathrm{full}}+\Delta t\,A_h^{\mathrm{full}})_{FF},\\
\mathbf b_h^i
&=
\bigl(
M_h^{\mathrm{full}}\mathbf p_{h,\mathrm{full}}^{i-1}
+\Delta t\,\mathbf f_{h,\mathrm{full}}^i
\bigr)_F\\
&\quad-
(M_h^{\mathrm{full}}+\Delta t\,A_h^{\mathrm{full}})_{FD}
\mathbf g_D.
\end{split}
\label{eq:ex2_discrete}
\end{equation}
Here, $\mathbf p_h^i$ contains the unknown pressures. For $i>1$, 
$\mathbf p_{h,\mathrm{full}}^{i-1}$ contains the
corrected solution from the preceding time step of the same
run, including its boundary values. The vector $\mathbf f_{h,\mathrm{full}}^i$
is the assembled load vector, and $\mathbf g_D$ contains
the prescribed pressure values.

We represent the coarse solution as $G_H+u_H^i$, where
$G_H(x,y)=1-x/L_x$ satisfies the prescribed pressure conditions
and $u_H^i$ vanishes on the Dirichlet boundaries.
We set $G_h=I_H^hG_H$ and let $\mathbf g_h$ contain its values
at the nodes indexed by $F$. With $P$ denoting the prolongation
matrix from coarse to fine unknowns, the Galerkin system is
\begin{equation}
B_H\mathbf u_H^i=\mathbf b_H^i,
\qquad
B_H=P^TB_hP,
\qquad
\mathbf b_H^i=P^T(\mathbf b_h^i-B_h\mathbf g_h).
\label{eq:ex2_galerkin}
\end{equation}
The matrix $B_H$ remains fixed throughout the simulation.

At each time step, the measured coarse solver produces the
reconstructed coefficients $\mathbf u_H^{q,i}$.
Let $u_H^{q,i}$ denote the corresponding finite element function.
Adding $G_H$ gives the reconstructed coarse solution
$\widetilde p_H^{q,i}=G_H+u_H^{q,i}$, including the prescribed
boundary values. The initial fine-grid vector is
\begin{equation}
\mathbf p_h^{(0),i}=\mathbf g_h+P\mathbf u_H^{q,i}.
\end{equation}
The classical solver (AMG-PCG) then solves
\begin{equation}
B_h\boldsymbol\delta_h^i
=\mathbf b_h^i-B_h\mathbf p_h^{(0),i}
\end{equation}
and gives the corrected solution
\begin{equation}
\mathbf p_h^{\mathrm{corr},i}
=\mathbf p_h^{(0),i}+\boldsymbol\delta_h^i.
\end{equation}
We use the same AMG hierarchy throughout each trajectory and
an explicitly recomputed relative residual tolerance of $10^{-10}$.

\paragraph{Measurement and parameter initialization}
We follow the measurement and reconstruction procedure of
Experiment~\ref{sec:ex1}. Each objective circuit uses 2,048
shots, and final tomography uses 4,096 shots for each of
the $3^4=81$ product-Pauli bases. All measurements are
simulated with Qiskit Aer. We use the same three measurement conditions as in Experiment~\ref{sec:ex1}:
finite-shot sampling without gate or readout noise, mild synthetic
noise, and device-modeled noise.

For each condition, we compare two strategies for initializing
the circuit parameters. Independent initialization uses fresh
parameters at every time step. With parameter reuse, we initialize the circuit parameters
at the first time step in the same way as for independent
initialization. At each subsequent step, we use the parameters
selected at the preceding step:
\begin{equation}
\boldsymbol{\theta}^{(0),i}
=\boldsymbol{\theta}^{*,i-1},
\qquad i>1,
\end{equation}
where $\boldsymbol{\theta}^{*,i-1}$ denotes the selected parameters
from the preceding step. We perform three independent runs of 20 time steps for each
combination of measurement condition and initialization
strategy, giving 18 runs in total.

For both strategies, we first reconstruct the coarse solution
using the initial circuit parameters. We accept these parameters without SPSA updates if the
relative residual of the reconstructed coarse solution satisfies
\begin{equation}
\frac{\|B_H\mathbf u_H^{q,i}-\mathbf b_H^i\|_2}
     {\|\mathbf b_H^i\|_2}
\leq 0.25.
\label{eq:ex2_coarse_acceptance}
\end{equation}
Otherwise, we perform SPSA updates, checking this criterion
every 20 updates, and stop when it is satisfied or the limit
of 100 updates is reached. We select
the tested parameters with the smallest measured coarse residual
and perform an independent final readout to reconstruct the
coarse solution used for fine-grid initialization.

We present the results in two parts. First, we assess
coarse-solution accuracy, state fidelity, and 
pressure errors using the nine runs with parameter
reuse. Second, we examine how parameter reuse affects sampling requirements and compare
AMG-PCG iterations for initialization from the reconstructed
coarse solution and the exact Galerkin solution.

\paragraph{i) Coarse-solution accuracy, state fidelity, and pressure errors}
We compare the initial and corrected fine-grid solutions with
an independently computed reference on the same spatial mesh
and time grid. Their relative $L^2$ errors are
\begin{equation}
E_{\mathrm{rel}}^{s,i}
=
\frac{\|p_h^{s,i}-p_h^{\mathrm{ref},i}\|_{L^2(\Omega)}}
     {\|p_h^{\mathrm{ref},i}\|_{L^2(\Omega)}},
\qquad
s\in\{(0),\mathrm{corr}\}.
\end{equation}
We evaluate these norms using the finite element mass matrix.

Figure~\ref{fig:ex2_time_pressure} shows the reconstructed
coarse solution, corrected fine-grid solution, and reference
solution at three times under device-modeled noise with
parameter reuse. We show the run with the median final
coarse residual among the three runs. The coarse solution
captures the high- and low-pressure regions but misses their
fine spatial features. Fine-grid correction recovers these
features and gives a solution close to the reference.

\begin{figure}[!h]
\centering
\includegraphics[width=\linewidth]{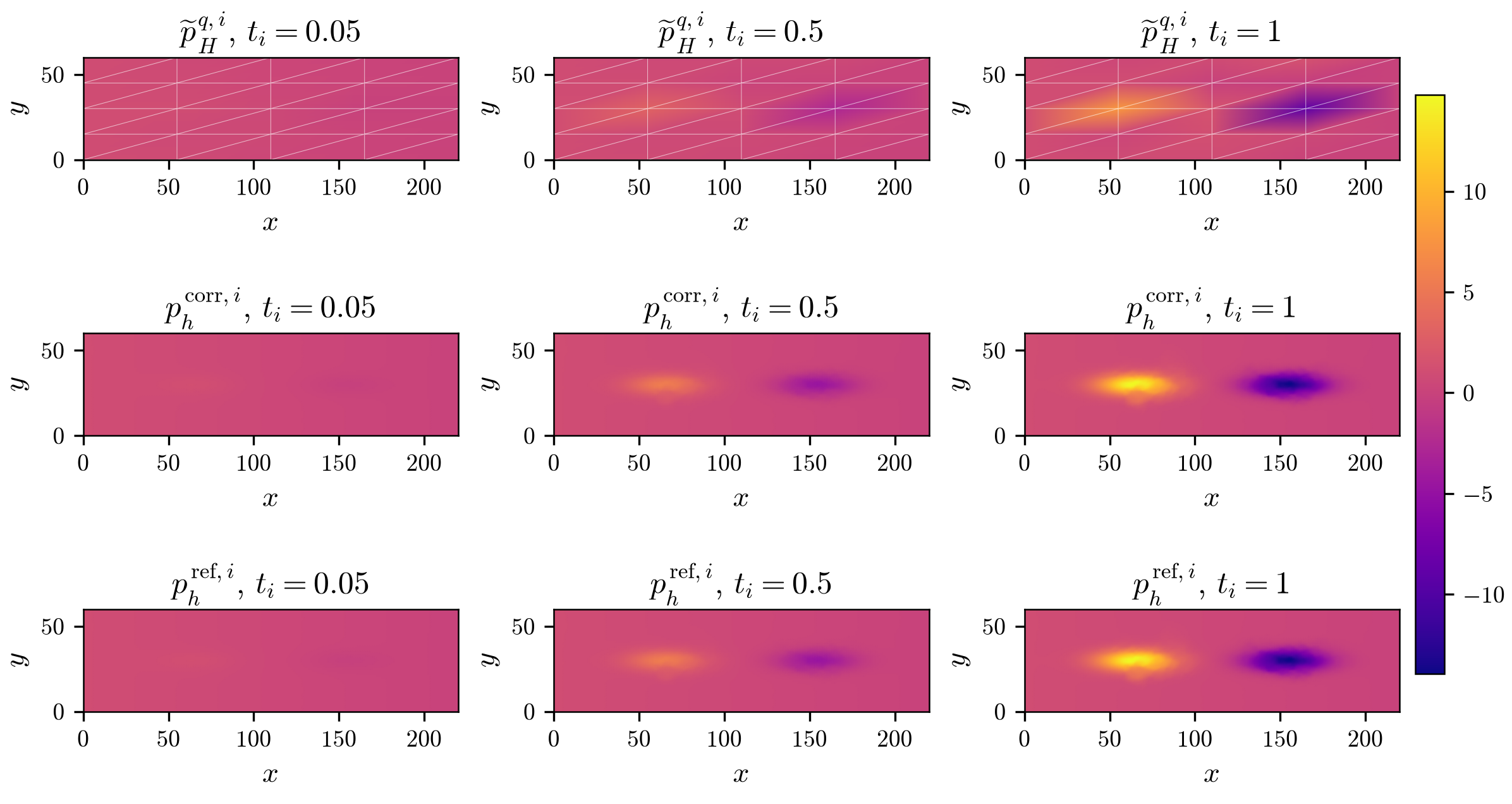}
\caption{Experiment 2. Pressure fields for the representative run under
device-modeled noise with parameter reuse. Rows show the reconstructed
coarse solution $\widetilde p_H^{q,i}$, corrected fine-grid
solution $p_h^{\mathrm{corr},i}$, and reference solution
$p_h^{\mathrm{ref},i}$. Columns correspond to
$t_i=0.05$, $0.5$, and $1$. All panels use the same color scale.
The coarse solution includes the prescribed boundary values,
and its mesh is shown in the top row.}
\label{fig:ex2_time_pressure}
\end{figure}

We next assess the reconstructed coarse solution and its
associated quantum state. The coefficient error measures
how closely the reconstructed coefficients match the exact
Galerkin solution. State fidelity measures how closely the
reconstructed quantum state matches the target state defined
by the normalized exact Galerkin solution.

To construct this target state, we compute
$\mathbf u_H^{\star,i}=B_H^{-1}\mathbf b_H^i$, append a zero
for the unused padded component, and normalize the resulting
vector. We define the coefficient error and state fidelity as

\begin{equation}
\begin{split}
E_{H,2}^i
&=
\frac{\|\mathbf u_H^{q,i}-\mathbf u_H^{\star,i}\|_2}
     {\|\mathbf u_H^{\star,i}\|_2},\\
|\psi_\star^i\rangle
&=
\frac{(\mathbf u_H^{\star,i},0)^T}
     {\|\mathbf u_H^{\star,i}\|_2},\\
\mathcal F_i
&=
\langle\psi_\star^i|\widehat\rho_i|\psi_\star^i\rangle.
\end{split}
\label{eq:ex2_quantum_metrics}
\end{equation}
Here, $\widehat\rho_i$ is the density matrix reconstructed
from independent final tomography measurements. A smaller
coefficient error indicates a more accurate coarse solution,
while fidelity $\mathcal F_i$ closer to one indicates better agreement
with the target state. The exact Galerkin solution is used
only to evaluate accuracy. These comparisons are made before
adding the boundary lift.

Figure~\ref{fig:ex2_coefficient_error} shows similar coefficient
errors under mild noise and device-modeled noise. Figure~\ref{fig:ex2_fidelity} shows that fidelity
remains near $0.89$ under finite-shot sampling and mild noise,
but near $0.61$ under device-modeled noise. Neither metric
changes noticeably at the source-amplitude change at $t=0.5$. 

\begin{figure}[!h]
\centering
\begin{subfigure}[t]{0.49\linewidth}
\centering
\includegraphics[width=\linewidth]
{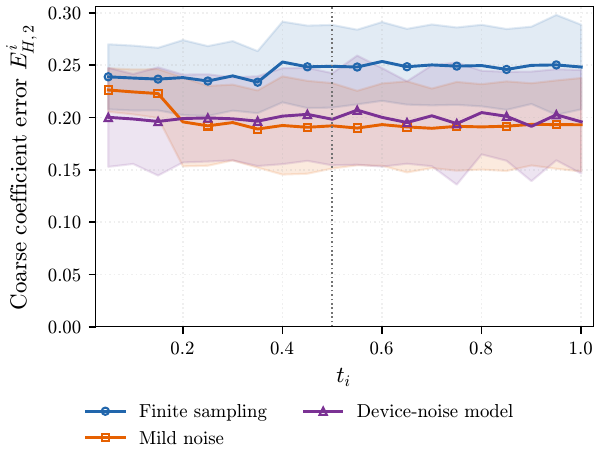}
\caption{Relative coefficient error $E_{H,2}^i$.}
\label{fig:ex2_coefficient_error}
\end{subfigure}
\hfill
\begin{subfigure}[t]{0.49\linewidth}
\centering
\includegraphics[width=\linewidth]
{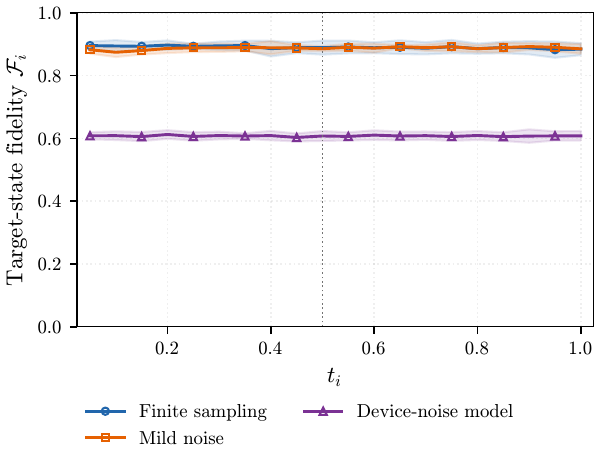}
\caption{Target-state fidelity $\mathcal F_i$.}
\label{fig:ex2_fidelity}
\end{subfigure}
\caption{Experiment 2. Coefficient error and state fidelity with parameter
reuse. Curves and shaded bands show means and sample standard
deviations over three trajectories per condition. The dotted
line marks the source-amplitude change at $t=0.5$.}
\label{fig:ex2_quantum_accuracy}
\end{figure}

Alongside these accuracy measures, we check whether probability
is assigned to the unused state introduced by padding.
The original coarse system has 15 coefficients, while its
four-qubit representation has 16 components. The extra
component, corresponding to $|1111\rangle$, should ideally
have zero probability. We estimate this probability as
\begin{equation}
\ell_{\mathrm{pad}}^i
=\frac{C_{1111}^i}{S_{\mathrm{pad}}^i},
\end{equation}
where $C_{1111}^i$ counts the $1111$ outcomes among the
$S_{\mathrm{pad}}^i$ shots of the final padding-measurement
circuit. With parameter reuse, the final-time estimates are
$(3.581\pm6.202)\times10^{-3}$ under finite-shot sampling,
$(1.628\pm1.973)\times10^{-3}$ under mild noise, and
$(4.069\pm0.282)\times10^{-3}$ under device-modeled noise.
These are means and sample standard deviations over three
runs per condition. Thus, at the final time, less than
$0.5\%$ of the measured probability is assigned to the unused
state on average under each condition.

Table~\ref{tab:ex2_quantum_pressure} compares pressure errors
before and after fine-grid correction at selected times.
The initial error increases over time, while the corrected
solution remains close to the discrete reference under
all three measurement conditions.

\begin{table}[!htbp]
\centering
\small
\setlength{\tabcolsep}{6pt}
\begin{tabular}{lrrr}
\toprule
Condition & $t_i$
& $E_{\mathrm{rel}}^{(0),i}$
& $E_{\mathrm{rel}}^{\mathrm{corr},i}$ \\
\midrule
Finite shots & 0.05
& $0.089\pm0.001$ & $1.49\times10^{-12}$ \\
& 0.5
& $0.465\pm0.005$ & $2.36\times10^{-11}$ \\
& 1
& $0.558\pm0.006$ & $3.90\times10^{-11}$ \\
\midrule
Mild noise & 0.05
& $0.089\pm0.001$ & $1.47\times10^{-12}$ \\
& 0.5
& $0.457\pm0.006$ & $2.44\times10^{-11}$ \\
& 1
& $0.549\pm0.007$ & $4.10\times10^{-11}$ \\
\midrule
Device-modeled noise & 0.05
& $0.088\pm0.001$ & $1.53\times10^{-12}$ \\
& 0.5
& $0.459\pm0.005$ & $2.45\times10^{-11}$ \\
& 1
& $0.550\pm0.008$ & $4.08\times10^{-11}$ \\
\bottomrule
\end{tabular}
\caption{Experiment 2. Relative $L^2$ pressure errors before and after
fine-grid correction with parameter reuse. Initial errors
are means $\pm$ sample standard deviations over three runs
per condition; corrected errors are maxima over those runs.}
\label{tab:ex2_quantum_pressure}
\end{table}

Thus, fine-grid correction gives comparable pressure accuracy
across the three measurement conditions, despite the lower
state fidelity under device-modeled noise. This accuracy
refers to agreement with the reference on the same spatial
mesh and time grid.


\paragraph{ii) Sampling requirements and AMG-PCG iterations}

We examine whether reusing circuit parameters reduces the sampling required for the coarse solve and whether coarse initialization reduces the iterations required for fine-grid correction. 

Table~\ref{tab:ex2_initialization_cost} compares the total
simulated shots per 20-step run for independent parameter
initialization and parameter reuse. The totals include gain
calibration, SPSA objective evaluations, intermediate
reconstructions, and independent final readouts.
These counts measure sampling requirements, not runtime.

\begin{table}[!h]
\centering
\small
\setlength{\tabcolsep}{5pt}
\begin{tabular}{llc}
\toprule
Condition & Parameter initialization
& Shots per run ($10^6$) \\
\midrule
Finite shots & Independent & $209.81\pm6.00$ \\
             & Parameter reuse & $25.03\pm4.33$ \\
Mild noise   & Independent & $212.41\pm3.00$ \\
             & Parameter reuse & $25.03\pm4.33$ \\
Device-modeled noise & Independent & $260.92\pm10.50$ \\
             & Parameter reuse & $26.54\pm2.60$ \\
\bottomrule
\end{tabular}
\caption{Experiment 2. Total simulated shots per 20-step run.
Entries are means $\pm$ sample standard deviations over
three trajectories per condition and initialization strategy.}
\label{tab:ex2_initialization_cost}
\end{table}

Parameter reuse reduces the mean sampling requirement by
approximately $88$-$90\%$ across the three measurement
conditions. Most time steps require no further SPSA updates:
56 of the 60 coarse solves under finite-shot sampling and
mild noise, and 57 under device-modeled noise. These solves
still require measurements to check the initial parameters
and reconstruct the coarse solution.

For fine-grid correction, the reconstructed coarse solution
requires 160 AMG-PCG iterations in total over 20 time steps
under all three measurement conditions and both parameter
strategies. Initialization from the exact Galerkin coarse
solution requires the same number of iterations. Both
initializations produce corrected solutions close to the
discrete reference. Thus, parameter reuse reduces sampling requirements.
Initialization from the reconstructed coarse solution gives
the same AMG-PCG iteration count and comparable corrected
pressure accuracy as exact Galerkin initialization.

\section{Conclusion}
\label{sec:conclusion}

{\color{black}
We developed a hybrid quantum-classical coarse-to-fine finite
element framework for heterogeneous Darcy flow. A variational
quantum solver approximates the Galerkin coarse solution;
measurement and classical reconstruction produce an explicit
coefficient vector. Prolongating that vector and adding the
boundary lifting gives an initial pressure for classical
AMG-PCG correction.
}

{\color{black}
The results in Experiment~1 show that measured quantum
initialization lowers the initial fine-grid $L^2$ and energy
errors relative to zero initialization under finite sampling
and modeled noise. However, this improvement does not reduce
the subsequent AMG-PCG iteration count. The fixed-parameter
study shows that additional shots improve coarse accuracy,
while modeled noise continues to affect the reconstruction.
In the time-dependent workflow, reusing circuit parameters
from the preceding time step reduces the mean sampling
requirements, while the total AMG-PCG iteration count
remains unchanged. This highlights the usefulness of
previous time-step information in reducing the measurement
cost of the quantum coarse solve.
}

{\color{black}
This framework provides a prototype for integrating quantum
approximations into established classical fine-grid solvers.
Although this work focuses on Darcy flow, the hybrid
coarse-to-fine PDE workflow can also be applied to other SPD
discretizations in multiphysics systems. Future work should
examine larger coarse systems, alternative coarse spaces,
more efficient reconstruction, and implementation on
quantum hardware.
}

\section*{Acknowledgments}
The work of all authors K. Chawla, S. Lee, and J. Valyou was supported by the U.S. Department of Energy, Office of Science, Energy Earthshots Initiatives under Award Number DE-SC 0024703.

\bibliographystyle{siamplain}
\bibliography{references} 

\end{document}

%% file: ex_shared.tex
\usepackage{amsfonts}
\usepackage{graphicx}
\usepackage{epstopdf}
\usepackage{algorithm}
\usepackage{algpseudocode} 
\newsiamremark{assumption}{Assumption}

\ifpdf
  \DeclareGraphicsExtensions{.eps,.pdf,.png,.jpg}
\else
  \DeclareGraphicsExtensions{.eps}
\fi

\newsiamremark{remark}{Remark}
\newsiamremark{hypothesis}{Hypothesis}
\crefname{hypothesis}{Hypothesis}{Hypotheses}
\newsiamthm{claim}{Claim}

\headers{Hybrid Quantum-Classical Coarse-to-Fine FEM Solver}{K. Chawla, S. Lee, and J. Valyou}

\title{A Hybrid Quantum-Classical Coarse-to-Fine Finite Element Solver for Heterogeneous Darcy Flow\thanks{Submitted to the editors.
\funding{This work was supported by the U.S. Department of Energy, Office of Science, Energy Earthshots Initiative under Award No. DE-SC0024703.}}}

\author{Kapil Chawla\thanks{Department of Mathematics, Florida State University, 1017 Academic Way, Tallahassee, FL 32306, United States of America (\email{kc25bc@fsu.edu}; \email{slee17@fsu.edu}; \email{jmv18b@fsu.edu}).}
\and Sanghyun Lee\footnotemark[2] \and Jonathan Valyou\footnotemark[2]}

\usepackage{amsopn}

\makeatletter
\newcommand*{\addFileDependency}[1]{
  \typeout{(#1)}
  \@addtofilelist{#1}
  \IfFileExists{#1}{}{\typeout{No file #1.}}
}
\makeatother


%% file: references.bib
@article{christie2001tenth,
  author  = {Christie, M. A. and Blunt, M. J.},
  title   = {Tenth {SPE} Comparative Solution Project: A Comparison of Upscaling Techniques},
  journal = {SPE Reservoir Evaluation \& Engineering},
  volume  = {4},
  number  = {04},
  pages   = {308--317},
  year    = {2001},
  doi     = {10.2118/72469-PA}
}

@Article{Nicolaides1987,
  title   = {Deflation of Conjugate Gradients with Applications to Boundary Value Problems},
  author  = {R. A. Nicolaides},
  journal = {SIAM Journal on Numerical Analysis},
  number  = {2},
  volume  = {24},
  year    = {1987},
  pages   = {355--365},
}

@Book{Saad2003,
  title     = {Iterative Methods for Sparse Linear Systems},
  author    = {Yousef Saad},
  publisher = {SIAM},
  address   = {Philadelphia, PA},
  edition   = {2},
  year      = {2003},
}

@Book{Smith1996,
  title     = {Domain Decomposition: Parallel Multilevel Methods for Elliptic Partial Differential Equations},
  author    = {Barry F. Smith and Petter E. Bj{\o}rstad and William D. Gropp},
  publisher = {Cambridge University Press},
  address   = {Cambridge, UK},
  year      = {1996},
}

@Book{ToselliWidlund2005,
  title     = {Domain Decomposition Methods---Algorithms and Theory},
  author    = {Andrea Toselli and Olof B. Widlund},
  publisher = {Springer},
  address   = {Berlin},
  series    = {Springer Series in Computational Mathematics},
  volume    = {34},
  year      = {2005},
}

@Article{Vuik2006,
  title   = {Deflation Techniques for Domain Decomposition Methods Applied to Porous Media Flow Problems},
  author  = {C. Vuik and R. Nabben and J. M. Tang},
  journal = {SIAM Journal on Scientific Computing},
  number  = {5},
  volume  = {27},
  year    = {2006},
  pages   = {1742--1766},
}

@Article{Harrow2009,
  title   = {Quantum Algorithm for Linear Systems of Equations},
  author  = {Aram W. Harrow and Avinatan Hassidim and Seth Lloyd},
  journal = {Physical Review Letters},
  number  = {15},
  volume  = {103},
  year    = {2009},
  pages   = {150502},
}

@Article{Kandala2017,
  title   = {Hardware-Efficient Variational Quantum Eigensolver for Small Molecules and Quantum Magnets},
  author = {Abhinav Kandala and Antonio Mezzacapo and Kristan Temme and Maika Takita and Markus Brink and Jerry M. Chow and Jay M. Gambetta},
  journal = {Nature},
  volume  = {549},
  year    = {2017},
  pages   = {242--246},
}

@InProceedings{Gilyen2019,
  title     = {Quantum Singular Value Transformation and Beyond: Exponential Improvements for Quantum Matrix Arithmetics},
  author    = {Andr\'{a}s Gily\'{e}n and Yuan Su and Guang Hao Low and Nathan Wiebe},
  booktitle = {Proceedings of the 51st Annual ACM SIGACT Symposium on Theory of Computing},
  year      = {2019},
  pages     = {193--204},
}

@Article{Lubasch2020,
  title   = {Variational Quantum Algorithms for Nonlinear Problems},
  author  = {Michael Lubasch and Jaewoo Joo and Pierre Moinier and Martin Kiffner},
  journal = {Physical Review A},
  number  = {1},
  volume  = {101},
  year    = {2020},
  pages   = {010301},
}

@Article{Liu2021,
  title   = {Variational Quantum Linear Solver for Poisson Equation},
  author  = {Hsin-Yuan Liu and Yu-Chen Wang and Yung-Cheng Luo and Ying-Jer Chen},
  journal = {Scientific Reports},
  volume  = {11},
  year    = {2021},
  pages   = {12453},
}

@Article{Balducci2022,
  title   = {Review and Perspectives in Quantum Computing for Partial Differential Equations in Structural Mechanics},
  author  = {Giorgio Tosti Balducci and Boyang Chen and Matthias M{\"o}ller and Marc Gerritsma and Roeland De Breuker},
  journal = {Frontiers in Mechanical Engineering},
  volume  = {8},
  year    = {2022},
}

@Article{Schalkers2024,
  title   = {Efficient and Fail-Safe Quantum Algorithm for the Transport Equation},
  author  = {Merel A. Schalkers and Matthias M{\"o}ller},
  journal = {Journal of Computational Physics},
  volume  = {502},
  year    = {2024},
  pages   = {112816},
}

@Article{Preskill2018,
  title   = {Quantum Computing in the {NISQ} Era and Beyond},
  author  = {John Preskill},
  journal = {Quantum},
  volume  = {2},
  year    = {2018},
  pages   = {79},
}

@Book{EfendievHou2009,
  title     = {Multiscale Finite Element Methods: Theory and Applications},
  author    = {Yalchin Efendiev and Thomas Y. Hou},
  publisher = {Springer},
  year      = {2009},
}

@Book{Briggs2000,
  title     = {A Multigrid Tutorial},
  author    = {William L. Briggs and Van Emden Henson and Steve F. McCormick},
  publisher = {SIAM},
  year      = {2000},
}

@Article{Fang2023,
  title   = {Qubit-Efficient Quantum Algorithm for Differential Equations},
  author  = {Yingkai Fang and others},
  journal = {Nature Communications},
  volume  = {14},
  year    = {2023},
  pages   = {1--9},
}

@Article{Durlofsky1991,
  title   = {Numerical Calculation of Equivalent Grid Block Permeability Tensors for Heterogeneous Porous Media},
  author  = {Louis J. Durlofsky},
  journal = {Water Resources Research},
  number  = {5},
  volume  = {27},
  year    = {1991},
  pages   = {699--708},
}

@Misc{Qiskit2021,
  title = {{Qiskit}: An Open-Source Framework for Quantum Computing},
  author = {{Qiskit contributors}},
  year = {2021},
}

@Misc{Morales2024QLSsurvey,
  title  = {Quantum Linear System Solvers: A Survey of Algorithms and Applications},
  author = {Mauro E. S. Morales and Lirand{\"e} Pira and Philipp Schleich and Kelvin Koor and Pedro C. S. Costa and Dong An and Al{\'a}n Aspuru-Guzik and Lin Lin and Patrick Rebentrost and Dominic W. Berry},
  year   = {2024},
  note   = {arXiv:2411.02522},
}

@Misc{Deiml2024,
  title  = {Quantum Realization of the Finite Element Method},
  author = {Matthias Deiml and Daniel Peterseim},
  year   = {2024},
  note   = {arXiv:2403.19512},
}

@Article{Arora2025,
  title   = {An Implementation of the Finite Element Method in Hybrid Classical/Quantum Computers},
  author  = {Abhishek Arora and Benjamin M. Ward and Caglar Oskay},
  journal = {Finite Elements in Analysis and Design},
  volume  = {248},
  year    = {2025},
  pages   = {104354},
}

@Misc{Jacob2025,
  title  = {Bridging Quantum and Classical Computing for Partial Differential Equations through Multifidelity Machine Learning},
  author = {Bruno Jacob and Amanda A. Howard and Panos Stinis},
  year   = {2025},
  note   = {arXiv:2512.05241},
}

@Article{LiZikatanovZuo2025,
  title   = {Reduced Krylov Basis Methods for Parametric Partial Differential Equations},
  author  = {Yuwen Li and Ludmil T. Zikatanov and Cheng Zuo},
  journal = {SIAM Journal on Numerical Analysis},
  number  = {2},
  volume  = {63},
  year    = {2025},
  pages   = {976--999},
}

@Article{BravoPrieto2023,
  title   = {Variational Quantum Linear Solver},
  author  = {Carlos Bravo-Prieto and Ryan LaRose and Marco Cerezo and Yi{\u g}it Suba{\c s}{\i} and Lukasz Cincio and Patrick J. Coles},
  journal = {Quantum},
  volume  = {7},
  year    = {2023},
  pages   = {1188},
}

@Misc{ChawlaLeeShin2026PAM,
title  = {A Parallel and Adaptive Mesh-Free Method for Discontinuous Coefficient Fields in Heterogeneous Porous Media},
author = {Kapil Chawla and Sanghyun Lee and Yeonjong Shin},
year   = {2026},
note   = {arXiv:2605.16564},
}

@article{mottonen2005,
  author  = {M{\"o}tt{\"o}nen, Mikko and Vartiainen, Juha J.
             and Bergholm, Ville and Salomaa, Martti M.},
  title   = {Transformation of quantum states using
             uniformly controlled rotations},
  journal = {Quantum Information and Computation},
  volume  = {5},
  number  = {6},
  pages   = {467--473},
  year    = {2005},
}

@article{Spall1992,
  author  = {Spall, James C.},
  title   = {Multivariate Stochastic Approximation Using
             a Simultaneous Perturbation Gradient Approximation},
  journal = {IEEE Transactions on Automatic Control},
  volume  = {37},
  number  = {3},
  pages   = {332--341},
  year    = {1992},
}

@article{James2001,
  author  = {James, Daniel F. V. and Kwiat, Paul G.
             and Munro, William J. and White, Andrew G.},
  title   = {Measurement of qubits},
  journal = {Physical Review A},
  volume  = {64},
  pages   = {052312},
  year    = {2001},
}

@article{Sato2021,
  author  = {Sato, Yuki and Kondo, Ruho and Koide, Satoshi
             and Takamatsu, Hideki and Imoto, Nobuyuki},
  title   = {Variational quantum algorithm based on the
             minimum potential energy for solving the
             Poisson equation},
  journal = {Physical Review A},
  volume  = {104},
  number  = {5},
  pages   = {052409},
  year    = {2021},
}

@article{Pollachini2021,
  author  = {Pollachini, Giovani G. and Salazar, Juan P. L. C.
             and Goes, Caio B. D. and Maciel, Thiago O.
             and Duzzioni, Eduardo I.},
  title   = {A hybrid classical-quantum approach to solve the
             heat equation using quantum annealers},
  journal = {Physical Review A},
  volume  = {104},
  number  = {3},
  pages   = {032426},
  year    = {2021},
}

@article{Aaronson2015,
  author  = {Aaronson, Scott},
  title   = {Read the fine print},
  journal = {Nature Physics},
  volume  = {11},
  number  = {4},
  pages   = {291--293},
  year    = {2015},
}

@article{Childs2017,
  author  = {Childs, Andrew M. and Kothari, Robin and Somma, Rolando D.},
  title   = {Quantum Algorithm for Systems of Linear Equations with
             Exponentially Improved Dependence on Precision},
  journal = {SIAM Journal on Computing},
  volume  = {46},
  number  = {6},
  pages   = {1920-1950},
  year    = {2017},

}

@article{Henderson2024,
  author  = {Henderson, Jessie M. and Kath, John and Golden, John K. and
             Percus, Allon G. and O'Malley, Daniel},
  title   = {Addressing quantum's ``fine print'' with efficient state
             preparation and information extraction for quantum algorithms
             and geologic fracture networks},
  journal = {Scientific Reports},
  volume  = {14},
  pages   = {3592},
  year    = {2024},
}

@book{ChenHuanMa2006,
  author    = {Chen, Zhangxin and Huan, Guanren and Ma, Yuanle},
  title     = {Computational Methods for Multiphase Flows in Porous Media},
  series    = {Computational Science and Engineering},
  volume    = {2},
  publisher = {Society for Industrial and Applied Mathematics},
  address   = {Philadelphia, PA},
  year      = {2006},
}

@article{Golden2022,
  author  = {Golden, John and O'Malley, Daniel and Viswanathan, Hari},
  title   = {Quantum computing and preconditioners for hydrological linear
             systems},
  journal = {Scientific Reports},
  volume  = {12},
  pages   = {22285},
  year    = {2022},
  
}

@article{Henderson2023,
  author  = {Henderson, Jessie M. and Podzorova, Marianna and Cerezo, M. and
             Golden, John K. and Gleyzer, Leonard and Viswanathan, Hari S. and
             O'Malley, Daniel},
  title   = {Quantum algorithms for geologic fracture networks},
  journal = {Scientific Reports},
  volume  = {13},
  pages   = {2906},
  year    = {2023},
  
}

@article{iten2016,
  author  = {Iten, Raban and Colbeck, Roger and Kukuljan, Ivan
             and Home, Jonathan and Christandl, Matthias},
  title   = {Quantum Circuits for Isometries},
  journal = {Physical Review A},
  volume  = {93},
  pages   = {032318},
  year    = {2016},
}

@misc{qiskitAerDepolarizing,
  author       = {{Qiskit Development Team}},
  title        = {Depolarizing Error},
  howpublished = {Qiskit Aer Documentation},
  note         = {Accessed October 16, 2026},
}

@misc{qiskitAerReadout,
  author       = {{Qiskit Development Team}},
  title        = {{ReadoutError}},
  howpublished = {Qiskit Aer Documentation},
  note         = {Accessed October 16, 2026},
}

@misc{qiskitFakeProvider,
  author       = {{IBM Quantum}},
  title        = {Fake Provider},
  howpublished = {Qiskit IBM Runtime Documentation},
  note         = {Accessed October 16, 2026},
}
